\documentstyle[amscd]{amsart}
\numberwithin{equation}{section}
\theoremstyle{plain}
\newtheorem{thm}{Theorem}[section]
\newtheorem{cor}[thm]{Corollary}
\newtheorem{lem}[thm]{Lemma}
\newtheorem{prop}[thm]{Proposition}
\begin{document}
\title{Full groups of one-sided topological Markov shifts}
\author{Kengo Matsumoto}
\address{ 
Department of Mathematics, 
Joetsu University of Education,
Joetsu 943-8512 Japan}
\email{kengo{@@}juen.ac.jp}

\begin{abstract}
Let $(X_A,\sigma_A)$
be the right one-sided topological Markov shift
for an irreducible matrix with entries in $\{0,1\}$,
and
$\Gamma_A$ the continuous full group of $(X_A,\sigma_A)$.
For two irreducible matrices $A$ and $B$ with entries in $\{0,1\}$,
it will be proved that the continuous full groups
$\Gamma_A$ and $\Gamma_B$ are isomorphic as abstract groups
 if and only if
their  one-sided topological Markov shifts
$(X_A,\sigma_A)$ and $(X_B,\sigma_B)$ are continuously orbit equivalent.   
\end{abstract}



\maketitle

%%%%%%%%%%%%%%%%%%%%%%%%%%%%%%%%%%%%%%%%%%
\def\Zp{{ {\Bbb Z}_+ }}
\def\M{{{\cal M}}}
\def\N{{{\cal M}}}
\def\H{{{\cal H}}}
\def\K{{{\cal K}}}
\def\A{{{\cal A}}}
\def\E{{ {\cal E} }}
\def\B{{{\cal B}}}
\def\R{{{\cal R}}}
\def\S{{{\cal S}}}
\def\Aut{{{\operatorname{Aut}}}}
\def\End{{{\operatorname{End}}}}
\def\Ext{{{\operatorname{Ext}}}}
\def\Hom{{{\operatorname{Hom}}}}
\def\Ker{{{\operatorname{Ker}}}}
\def\Index{{{\operatorname{Index}}}}
\def\index{{{\operatorname{index}}}}
\def\id{{{\operatorname{id}}}}
\def\Ad{{{\operatorname{Ad}}}}
\def\dim{{{\operatorname{dim}}}}
\def\Homeo{{{\operatorname{Homeo}}}}
\def\supp{{{supp}}}
\def\Proj{{{Proj}}}
\def\span{{{span}}}
\def\ProbA{{{Prob(X_A)}}}
\def\ProbB{{{Prob(X_B)}}}
\def\OA{{{\cal O}_A}}
\def\OB{{{\cal O}_B}}
\def\DA{{{\cal D}_A}}
\def\DB{{{\cal D}_B}}
\def\FA{{{\cal F}_A}}
\def\FB{{{\cal F}_B}}
%%%%%%%%%%%%%%%%%%%%%%%%%%%%%%%%%%%%%%% 

 %%%%%%%%%%%%%%%%%%%%%%%%%%%%%%%%%%%%%%%%%%%%%%%%%%%%%%%%% 
\section{Introduction}
%%%%%%%%%%%%%%%%%%%%%%%%%%%%%%%%%%%

Giordano-Putnam-Skau have introduced studies of orbit equivalences for 
minimal homeomorphisms on Cantor sets (\cite{GPS}, \cite{GPS2}, 
cf. \cite{GMPS}, \cite{HPS}, \cite{Put}, etc.).
Minimal homeomorphisms on Cantor sets are now called Cantor minimal systems.
In their theory, the full groups   
play crucial r\^{o}le to classify Cantor minimal systems under orbit equivalences.
The full groups are defined as groups of homeomorphisms of the Cantor sets 
whose orbits are contained in the orbits of the minimal homeomorphisms.
%The classification of Cantor minimal systems are closely reletaed to
%that of the associated $C^*$-crossed products.
Giordano-Putnam-Skau have proved thet the full groups as groups 
 are complete invariants for orbit equivalences of Cantor minimal systems
(\cite{GPS2}). 
(For measure theoretic studies for orbit equivalences of ergodic transformations,
see \cite{CoKr},  \cite{D}, \cite{D2}, \cite{HO}, \cite{Kr}, etc.).

The class of topological Markov shifts
is another important class
of topological dynamical systems on Cantor sets. 
The author has introduced a study of orbit equivalence of one-sided topological Markov shifts, related to classification of Cuntz-Krieger algebras 
(\cite{MaPacific}, \cite{MaPAMS}).  
Let $A$ be an $N\times N$ irreducible matrix with entries in $\{0,1\}$
satisfying  condition (I) in the sense of Cuntz-Krieger \cite{CK}.
Let us denote by $X_A$ the shift space
$$
X_A = \{ (x_n )_{n \in \Bbb N} \in \{1,\dots,N \}^{\Bbb N}
\mid
A(x_n,x_{n+1}) =1 \text{ for all } n \in {\Bbb N}
\}
$$
over $\{1,\dots,N\}$.
It is homeomorphic to a Cantor set in  natural  product topology.
The continuous surjective map $\sigma_A$ on $X_A$ is defined by 
$\sigma_{A}((x_n)_{n \in {\Bbb N}})=(x_{n+1} )_{n \in {\Bbb N}}$.
The topological dynamical system 
$(X_A, \sigma_A)$ is called the (right) one-sided topological Markov shift for $A$.
The continuous map $\sigma_A$ on $X_A$ 
is no longer homeomorphism
and has dense periodic points.
Hence the one-sided topological Markov shifts are considered to
locate far from Cantor minimal systems
in the topological dynamical systems on Cantor sets.
The continuous full group $\Gamma_A$ of
$(X_A,\sigma_A)$
is defined to be the group of homeomorphisms
$\tau$ on $X_A$ satisfying 
$$
\tau(x) \in orb_{\sigma_A}(x) = 
\cup_{k=0}^\infty \cup_{l=0}^\infty \sigma_A^{-k}(\sigma_A^l(x))
\quad\text{ 
for all } x \in X_A
$$
and having continuous orbit cocycles.
%The group $\Gamma_A$ naturally acts on $X_A$. 

In this paper, 
we will study the continuous full groups on one-sided topological Markov shifts
by using techniques in \cite{GPS2} for analyzing full groups of Cantor minimal systems.  
Keeping in mind the results of the Giordano-Putnam-Skau's paper \cite{GPS2} on the Cantor minimal systems,
the algebraic structure of the continuous full groups 
are expected to determine the structure of continuous orbit equivalences of one-sided topological Markov shifts.
Let $A,B$ be two irreducible square matrices 
with entries in $\{0,1\}$ satisfying condition (I).
If there exists a homeomorphism
$h: X_A \longrightarrow X_B$ 
 such that
$
h \circ \Gamma_A \circ h^{-1} = \Gamma_B,
$
then the continuous full groups $\Gamma_A$ 
and
 $\Gamma_B$ on $X_B$
are said to be  spatially isomorphic. 
We have proved in \cite{MaPacific} that 
$\Gamma_A$ and $\Gamma_B$ are spatially isomorphic 
if and only if the one-sided topological Markov shifts
$(X_A,\sigma_A)$ and $(X_B,\sigma_B)$ are 
continuously orbit equivalent.
We will prove in this paper 
that 
an algebraic isomorphism between the continuous full groups of 
one-sided topological Markov shifts  yields an spatial isomorphism between them.
Hence the algebraic structure of the continuous full groups of one-sided topological Markov shifts
 are complete invariants of 
continuous orbit equivalence classes of one-sided topological Markov shifts.
This is 
 an anologue of the Giordano-Putnam-Skau's theorem \cite[Theorem 4.2]{GPS2}
for one-sided topological Markov shifts.
The sterategy to prove it basically follows Giordano-Putnam-Skau's paper \cite{GPS2} 
in which
the algebraic isomorphism between  the full groups of Cantor minimal systems 
induces a spatial isomorphism between them
(\cite[Theorem 4.2]{GPS2}).
The continuous full groups of one-sided topological Markov shifts 
are countable, nonamenable groups (cf. \cite{MaPre2012}).
They are huge groups 
rather than the full groups of Cantor minimal systems. 
Indeed the full groups of Cantor minimal systems have invariant probability measures,
whereas the continuous full groups of one-sided topological Markov shifts
do not have any invariant probability measure.
In following the proof of \cite[Theorem 4.2]{GPS2}, 
there are several places where the proofs of \cite[Theorem 4.2]{GPS2} do not well work 
in our setting because our full groups are huge.      
We may modify their proofs which we apply to our situations and      
will reach our goal stated as the following theorem:
\begin{thm} \label{thm:main}
Let $A,B$ be two irreducible square matrices 
with entries in $\{0,1\}$ satisfying condition (I).
Then the following two conditions are equivalent:
\begin{enumerate} 
\renewcommand{\labelenumi}{(\roman{enumi})}
\item
The continuous full groups $\Gamma_A$ and $\Gamma_B$ are isomorphic
as abstract groups.
\item
The continuous full groups $\Gamma_A$ and $\Gamma_B$ are spatially isomorphic.
\end{enumerate}
\end{thm}
As (ii) $\Longrightarrow$ (i) is clear,
the implication  (i) $\Longrightarrow$ (ii) is the main part. 
For an irreducible square matrix $A$ with entries in $\{0,1\}$,
denote by
$\OA$ its Cuntz-Krieger algebra.
We also detnote by
$\DA$ its canonical maximal abelian subalgbra of $\OA$.   
The class of Cuntz-Krieger algebras plays an important r\^{o}le in classification theory of simple purely infinite $C^*$-algebras. 
Related to the classification of the Cuntz-Krieger
algebras, the author has shown that
there exists  an isomorphism $\Psi: \OA \rightarrow  \OB$
 satisfying $\Psi(\DA) = \DB$
if and only if 
 $\Gamma_A$ 
and $\Gamma_B$ are spatially isomorphic,
which is also equivalent to the condition that 
the one-sided topological Markov shifts
$(X_A,\sigma_A)$ and $(X_B,\sigma_B)$  
are continuously orbit equivalent (\cite{MaPacific}).
Therefore we know
\begin{cor}
Let $A,B$ be two irreducible square matrices 
with entries in $\{0,1\}$ satisfying condition (I).
Then the following three conditions are equivalent:
\begin{enumerate} 
\renewcommand{\labelenumi}{(\roman{enumi})}
\item
The continuous full groups $\Gamma_A$ 
and $\Gamma_B$  are isomorphic
as abstract groups.
\item
The one-sided topological Markov shifts
$(X_A,\sigma_A)$ and $(X_B,\sigma_B)$  
are continuously orbit equivalent.
\item
There exists an isomorphism $\Psi: \OA \rightarrow  \OB$
 satisfying $\Psi(\DA) = \DB$.
\end{enumerate}
\end{cor}
Let $N$ and $M$ be the size of matrix $A$ and that of $B$ respectively.
Denote by $I_N$ and by $I_M$ the identity matrix of size $N$ and that of size $M$ respectively.  
In \cite{MaPAMS},
under the condition that
$\det(A- I_N)\det(B-I_M) \ge 0$,
an isomorphism between Cuntz-Krieger algebras 
induces an isomorphism between them which preserves their canonical 
maximal abelian subalgebras.
Hence we have
\begin{cor}
Let $A,B$ be two irreducible square matrices 
with entries in $\{0,1\}$ satisfying condition (I).
Suppose that
$\det(A-I_N)\det(B- I_M) \ge 0$.
Then the continuous full  groups $\Gamma_A$ and $\Gamma_B$ 
 are isomorphic as abstract groups
if and only if the Cuntz-Krieger algebras
$\OA$ and $\OB$ are isomorphic.
\end{cor}
By using classification theorem for Cuntz-Krieger algebras 
obtained by M. R{\o}rdam
\cite{Ro}(cf. \cite{Ro2}),
one may classify the continuous full groups 
in terms of the underlying matrices under the determinant contition
$\det(A-I_N) \det(B-I_M) \ge 0$ as follows:
\begin{cor}
Let $A,B$ be two irreducible square matrices 
with entries in $\{0,1\}$ satisfying condition (I).
Suppose that
$\det(A- I_N)\det(B- I_M) \ge 0$.
The continuous full groups $\Gamma_A$ and $\Gamma_B$ 
 are isomorphic as abstract groups
 if and only if
 there exists an isomorphism
$\varPhi: {\Bbb Z}^N / (A^t -I_N) {\Bbb Z}^N
 \longrightarrow
  {\Bbb Z}^M / (B^t -I_M) {\Bbb Z}^M
  $
  such that
$\varPhi([1,\dots,1]) = [1,\dots,1].$
\end{cor}
Therefore we know that there are many mutually nonisomorphic 
continuous full groups of one-sided topological Markov shifts
such as the following corollary. 
\begin{cor}
Let $N, M $ be positive integers such that $N,M > 1$.
Denote by 
$\Gamma_{[N]}$ and $ \Gamma_{[M]}$ the continuous full groups of the one-sided full $N$-shift and $M$-shift respectively.  
Then 
$\Gamma_{[N]}$ and $ \Gamma_{[M]}$ are isomorphic as abstract groups if and only if
$N = M$.
\end{cor}
The main part of the paper is devoted to proving the implication 
(i) $\Longrightarrow$ (ii)
of Theorem \ref{thm:main}.
The paper is organized to prove it 
as in the following way.
In Section 2,  some basic properties of continuous full groups will be stated.
In Section 3, an open set of $X_A$ will be described in termes of a pair, called a strong commuting pair,  of subgroups of $\Gamma_A$.
In Section 4, a condition under which 
an open set of $X_A$ becomes clopen 
will be described in terms of algebraic conditions of a strong commuting pair. 
In Section 5, support of a strong commuting pair will be defined and  proved to be clopen.
In Section 6, an clopen set of $X_A$ 
will be completely replaced in terms of a pair of subgroups of $\Gamma_A$
with some algebraic conditions.
The pair of subgroups are called Dye pairs.
In Section 7, the main result and its corollaries will be stated.
The above sterategy to prove Theorem \ref{thm:main} basically  follows \cite{GPS2}.
However several proofs of in particular Lemma \ref{lem:3.18}, Lemma \ref {lem:3.19}, Lemma \ref {lem:3.21} and 
Lemma \ref{lem:3.24} are essentially different from the paper \cite{GPS2}.
The discussions of Section 4 and Section 6 are tough parts in this paper.

%%%%%%%%%%%%%%%%%%%%%%%%%%%%
%%%%%%%%%%%%%%%%%%%%%%%%%%%%%%%%%%%%%%%%%%%
\section{The continuous full groups and the local subgroups}
%%%%%%%%%%%%%%%%%%%%%%%%%%%%%%%%%%%%%%%%%%%%%

Let $A=[A(i,j)]_{i,j=1}^N$ 
be an $N\times N$ matrix with entries in $\{0,1\}$,
where $1< N \in {\Bbb N}$.
The matrix $A$ is always assumed to be essential,
which means that it has no zero columns or zero rows.
We assume that 
$A$ is irreducible and 
satisfies condition (I) in the sense of Cuntz-Krieger \cite{CK}.
In what follows, we fix the matrix $A$.
We denote by 
$X_A$ the shift space 
$$
X_A = \{ (x_n )_{n \in \Bbb N} \in \{1,\dots,N \}^{\Bbb N}
\mid
A(x_n,x_{n+1}) =1 \text{ for all } n \in {\Bbb N}
\}
$$
over $\{1,\dots,N\}$
of the right one-sided topological Markov shift for $A$.
It is a compact Hausdorff space in natural  product topology.
The condition (I) for $A$
is equivalent to the condition that $X_A$ is homeomorphic to a Cantor discontinuum.
The shift transformation $\sigma_A$ on $X_A$ is defined by 
$\sigma_{A}((x_n)_{n \in {\Bbb N}})=(x_{n+1} )_{n \in {\Bbb N}}$
for $(x_n)_{n \in {\Bbb N}}.$
It is a continuous surjective map on $X_A$.
The topological dynamical system 
$(X_A, \sigma_A)$ is called the (right) one-sided topological Markov shift for $A$.
%Throughout the paper, 
%the matrix $A$ is assumed to be irreducible and satisfy condition (I). 
A word $\mu = (\mu_1, \dots, \mu_k)$ for $\mu_i \in \{1,\dots,N\}$
is said to be admissible in $X_A$ 
if $\mu$ appears in somewhere in some element $x$ in $X_A$.
The length of $\mu$ is $k$, which is denoted by $|\mu|$.
 We denote by 
$B_k(X_A)$ the set of all admissible words of length $k \in {\Bbb N}$.
For $k=0$ we denote by $B_0(X_A)$ the empty word $\emptyset$.
We set 
$B_*(X_A) = \cup_{k=0}^\infty B_k(X_A)$ 
the set of admissible words of $X_A$.
For $x = (x_n )_{n \in {\Bbb N}} \in X_A$ 
and positive integers $k,l$ with
$k \le l$, 
we put the word 
$x_{[k,l]} = (x_k,x_{k+1},\dots, x_{l}) \in B_{l-k+1}(X_A)$
and the right infinite sequence
$x_{[k,\infty)} =(x_k, x_{k+1}, \dots ) \in X_A$.
We similarly use the notation 
$\mu_{[k,l]} = (\mu_k, \mu_{k+1},\dots, \mu_{l}) \in B_{l-k+1}(X_A)$ 
for a word $\mu= (\mu_1,\dots,\mu_m)\in B_m(X_A)$ with $k\le l\le m$.
For words 
$\mu = (\mu_1,\cdots,\mu_k) \in B_k(X_A),
\nu = (\nu_1,\cdots,\nu_l) \in B_l(X_A)$
with $A(\mu_k,\nu_1) =1$,
denote by $\mu\nu$
its concatenation
$
\mu\nu = (\mu_1,\cdots,\mu_k,\nu_1,\cdots,\nu_l) \in B_{k+l}(X_A).
$
For a word
$\mu = (\mu_1,\cdots,\mu_k) \in B_k(X_A)$,
we denote by $U_\mu$ its cylinder set
$$
U_\mu = \{(x_n)_{n\in{\Bbb N}} \in X_A \mid x_1=\mu_1,\dots,\mu_k = x_k \}.
$$

For $x = (x_n )_{n \in \Bbb N} \in X_A$,
the orbit $orb_{\sigma_A}(x)$ of $x$ under $\sigma_A$ 
is defined by
$$
orb_{\sigma_A}(x) = 
\cup_{k=0}^\infty \cup_{l=0}^\infty \sigma_A^{-k}(\sigma_A^l(x)) \subset X_A.
$$
Hence  
$ y =( y_n )_{n \in \Bbb N} \in X_A$ 
belongs to $orb_{\sigma_A}(x)$ 
if and only if
there exist
$k,l \in \Zp$ and an admissible word 
$(\mu_1, \dots, \mu_k) \in B_k(X_A)$ 
such that 
$$
y = (\mu_1,\dots, \mu_k, x_{l+1}, x_{l+2},\dots ).
$$
We denote by $\Homeo(X_A)$
the group of all homeomorphisms on $X_A$.
We have defined 
in \cite{MaPacific}
 the continuous  full group
$\Gamma_A$
for $(X_A,\sigma_A)$
as in the following way.

\noindent
{\bf Definition }(\cite{MaPacific}).
Let
$
\tau 
$ 
be a homeomorphism
on
$X_A$ 
such that
$
\tau(x) \in  orb_{\sigma_A}(x)
$ 
for all
$ 
x \in X_A.
$ 
Hence there exist functions
$k, l : X_A \rightarrow \Zp$ 
such that 
\begin{equation}
\sigma_A^{k(x)}(\tau(x) )=\sigma_A^{l(x)}(x)
\quad\text{ 
for all } x \in X_A. \label{eqn:cfg}
\end{equation}
Let $\Gamma_A$ 
be the set of all homeomorphisms $\tau$
such that  
there exist continuous functions 
$k, l : X_A \rightarrow \Zp$ 
satisfying \eqref{eqn:cfg}.
The set $\Gamma_A$ is a subgroup of $\Homeo(X_A)$ and is called 
the {\it continuous full group} of  $(X_A,\sigma_A)$.
The functions $k,l$ above are called the orbit cocycles for $\tau$.
They are not necessarily uniquely determined  
by $\tau$.
We note that the group
$\Gamma_A$ has been written as $[\sigma_A]_c$ in the earlier papers
\cite{MaPacific}, \cite{MaPAMS}.

A continuous map $\tau:X_A \longrightarrow X_A$
is called a {\it cylinder map}\,
if there exist $L \in {\Bbb N}$ 
and $\Phi : B_L(X_A) \longrightarrow B_*(X_A)$
such that
$\Phi(\nu) = \mu(\nu)= (\mu_1(\nu),\dots,\mu_{k(\nu)}(\nu)) \in B_{k(\nu)}(X_A)$
for 
$\nu = \nu_1\cdots\nu_L 
\in B_L(X_A)$
satisfies
%$\tau: U_\nu \longrightarrow U_{\Phi(\nu)}$ satisfies
\begin{equation*}
\tau(\nu_1,\dots, \nu_L, x_{L+1}, x_{L+2}, \dots ) 
=
(\mu_1(\nu),\dots,\mu_{k(\nu)}(\nu),  x_{L+1}, x_{L+2}, \dots ) 
\label{eqn:cylinder}
\end{equation*}
for
$( x_{L+1}, x_{L+2}, \dots ) \in X_A$
with
$A(\nu_L, x_{L+1}) =1$.
That is,
$\tau$ satisfies
% $$ \tau:  X_A = \sqcup_{\nu \in B_L(X_A)} U_\nu 
%\longrightarrow 
%\sqcup_{\nu \in B_L(X_A)} U_{\Phi(\nu)}
%\text{ such that }
$
\tau(U_\nu) = U_{\Phi(\nu)}
$
for $\nu \in B_L(X_A)$
with a map
$\Phi : B_L(X_A) \longrightarrow B_*(X_A)$.
We see that a homeomorphism $\tau$ of $X_A$
belongs to
$\Gamma_A$ if and only if 
$\tau$ is a cylinder map
(\cite{MaPre2012}).

We will first study
local structure of elements of  $\Gamma_A$.
Following \cite{GPS2}, we will use the notations as follows:
\begin{align*}
O(X_A) & = \text{The set of all open sets of } X_A. \\
CL(X_A) & = \text{The set of all closed sets of } X_A. \\
CO(X_A) & = \text{The set of all clopen sets of } X_A. 
\end{align*}
We note that an open set of $X_A$ is a countable disjoint union of cylinder sets
and a clopen set of $X_A$ is a finite disjoint union of cylinder sets.
\begin{lem}\label{lem:2.1}
For nonempty open sets
$U, Y \in O(X_A)$ and $x \in U$,
there exist $V\in CO(X_A)$ and
$\alpha \in \Gamma_A$
such that
\begin{equation*}
x \in V \subset U, \quad
\alpha(V) \subset Y, \quad
\alpha^2 = \id, \quad
\alpha|_{(V \cup\alpha(V))^c} =\id.
\end{equation*}
\end{lem}
\begin{pf}
Take a cylinder set
$U_\mu$ for some word $\mu = (\mu_1,\dots, \mu_n) \in B_n(X_A)$
such that $U_\mu \subset Y$.
Since $A$ satisfies condition (I),
there exist distinct words
$s=(s_1,\dots,s_k),s' = (s'_1,\dots,s'_k) \in B_k(X_A)$
and $u \in \{1,\dots,N\}$ such that
$$
A(\mu_n,s_1) =A(\mu_n,s'_1) =A(s_k,u) =A(s'_k,u) =1.
$$
There exists a word $\nu =(\nu_1,\dots,\nu_m) \in B_m(X_A)$
such that
$m > n + k + 1$ and $x \in U_\nu \subset U$.
Now $A$ is irreducible so that
one may find
a word
$\xi = (\xi_1,\dots,\xi_l) \in B_l(X_A)$
such that
$A(u,\xi_1) = A(\xi_l,\nu_m) =1$.
Put 
\begin{align*}
\bar{\mu} & 
= (\mu_1,\cdots, \mu_n, s_1,\dots,s_k, u, \xi_1,\cdots,\xi_l, \nu_m) 
\in B_{n + k+l+2}(X_A),\\
\bar{\mu}^\prime & 
= (\mu_1,\cdots, \mu_n, s'_1,\dots,s'_k, u, \xi_1,\cdots,\xi_l, \nu_m) 
\in B_{n + k+l+2}(X_A).
\end{align*}
Since
$\bar{\mu}_{[1,n+k+1]} \ne \bar{\mu}^\prime_{[1,n+k+1]}$
and
$|\nu |=m > n+k+1$,
at least either 
$\bar{\mu}_{[1,n+k+1]}$
or
$\bar{\mu}^\prime_{[1,n+k+1]}$
is different from
$\nu_{[1,n+k+1]}$.
We assume that
$\bar{\mu}_{[1,n+k+1]}
\ne \nu_{[1,n+k+1]}$,
so that
$U_\nu \cap U_{\bar{\mu}} = \emptyset$.
Put
$V = U_\nu$ and 
$L = n + k+ l +2$.
Define
$\alpha \in \Gamma_A$ by setting
\begin{equation*}
\alpha(x)
=
\begin{cases}
\bar{\mu} x_{[m+1,\infty)} 
& \text{ if } x = \nu x_{[m+1,\infty)} \in U_\nu, \\
\nu x_{[L+1,\infty)} 
& \text{ if } x = \bar{\mu} x_{[L+1,\infty)} \in U_{\bar{\mu}}, \\
x 
& \text{ otherwise} 
\end{cases}
\end{equation*}
for $x \in X_A$.
Then we have
\begin{equation*}
\alpha(V) =  U_{\bar{\mu}} \subset U_\mu \subset Y,
\qquad
\alpha^2 = \id, \qquad
\alpha|_{(U_\nu \cup U_{\bar{\mu}})^c} = \id. 
\end{equation*}
\end{pf}
For clopen sets $U,V$ of $X_A$, 
if there exists $\gamma \in \Gamma_A$ 
such that 
$\gamma(U) = V$,
then $U$ is said to be $\Gamma_A$-equivalent to $V$
and
written as
$U \underset{\Gamma_A}{\sim} V$.
\begin{lem}\label{lem:2.2}
Let 
$U, V \in CO(X_A)$
be nonempty clopen sets such that $U \cap V =\emptyset$.
If $U \underset{\Gamma_A}{\sim} V$,
there exists 
$\alpha \in \Gamma_A$
such that
\begin{equation}
\alpha(U) = V, \quad
\alpha^2 = \id, \quad
\alpha|_{(U \cup V)^c} =\id. \label{eqn:2.2}
\end{equation}
\end{lem}
\begin{pf}
Let $\gamma \in \Gamma_A$ satisfy
$\gamma(U) = V$.
One may define 
$\alpha \in \Gamma_A$ 
by setting
\begin{equation*}
\alpha(x) =
\begin{cases}
\gamma(x) & \text{ if } x \in U, \\
\gamma^{-1}(x) & \text{ if } x \in V, \\
x & \text{ otherwise}
\end{cases}
\end{equation*}
for $x \in X_A$.
As both $U$ and $V$ are clopen, 
$\alpha$ defines an element of $\Gamma_A$
which satisfies \eqref{eqn:2.2}.
\end{pf}
\begin{lem}[cf. {\cite[Lemma 3.4]{GPS2}}] \label{lem:3.4}
For any 
$U \in CO(X_A)$
and $x \in U$,
there exists 
$\alpha \in \Gamma_A$
such that
\begin{equation*}
\alpha(x) \ne x, \quad
\alpha^2 = \id, \quad
\alpha|_{U^c} =\id.
\end{equation*}
\end{lem}
\begin{pf}
Take clopen sets 
$U_1, Y_1 \subset U$
such that
$x \in U_1$ and $U_1 \cap Y_1 = \emptyset.$
By
Lemma \ref{lem:2.1},
there exist
a clopen set
$V_1$ with $x \in V_1 \subset U_1$
and 
$\alpha \in \Gamma_A$
such that
\begin{equation*}
\alpha(V_1) \subset Y_1, \quad
\alpha^2 = \id, \quad
\alpha|_{(V_1\cup \alpha(V_1))^c} =\id.
\end{equation*}
Since
$V_1 \cap Y_1 =\emptyset$
and
$U_1 \cup Y_1 \subset U$,
one has 
$\alpha(x) \ne x$
and
$\alpha|_{U^c} =\id.$
\end{pf}
We follow the notations below from \cite{GPS2}.

\noindent
{\bf Definition.}
For an open set $O$ of $X_A$, we set
\begin{equation*}
\Gamma_O = \{ \gamma \in \Gamma_A \mid \gamma(x) = x \text{ for all }
x \in O^c \}.
\end{equation*}
 A subgroup of $\Gamma_A$ of the form
$\Gamma_U$ for $U \in CO(X_A)$ is called a {\it local subgroup } of $\Gamma_A$.
We note that the notation $\Gamma_A$
for matrix $A$ is fixed and always denoting the continuous full group of $(X_A,\sigma_A)$. 
For a subgroup $H$ of $\Gamma_A$ the commutant of $H$ will be denoted by
$H^{\perp}$:
\begin{equation*}
H^{\perp} = \{ \xi \in \Gamma_A \mid \xi \gamma = \gamma \xi
\text{ for all } \gamma \in H \}
\end{equation*}
which is a subgroup of $\Gamma_A$.
The following proposition shows that
the action of local subgroups
$\Gamma_O$ on the underlying space $X_A$ 
is different from that of loclal subgroups of the full groups of Cantor minimal systems.
\begin{prop} \label{prop:2.4}
Let $O$ be a nonempty open set of $X_A$.
There is no regular Borel probability measure on $O$ invariant under $\Gamma_O$. 
\end{prop}
\begin{pf}
Suppose that there is a regular Borel probability measure 
$\mu$ on $O$ such that
$\mu \circ g = \mu$ for all $g \in \Gamma_O$.
Take $\nu(i) \in B_*(X_A), i=1,2,\dots$ 
such that
$O$ is a disjoint union $\sqcup_{i=1}^\infty U_{\nu(i)}$ so that
$
1 = \mu(O) = \sum_{i=1}^\infty \mu(U_{\nu(i)}).
$
One may find a word $\nu \in B_*(X_A)$
such that $U_\nu \subset O$ and
$\mu(U_\nu) >0$.
Take $x=(x_n)_{n\in {\Bbb N}} \in U_\nu$
such that 
$\mu(U_{x_{[1,n]}}) >0$ for all $n \in \Bbb N$
where
$x_{[1,n]} =(x_1,\dots, x_n)$.
One may find
$i >|\nu|$ and $k >1$ such that
$x_i = x_{i+k}$. 
Put $u = x_i$ and
\begin{align*} 
\zeta &  = (x_1,\dots, x_{i-1}), \quad \qquad \,
\xi  = (x_{i+1}, \dots, x_{i+k-1}), \\
\bar{\zeta} &  = (x_1,\dots, x_{i-1}, x_i), \qquad
\bar{\xi} = (x_{i+1}, \dots, x_{i+k-1}, x_{i+k})
\end{align*}
so that
$\bar{\zeta} \bar{\xi} = x_{[1,i+k]}$.
As
$i-1 \ge |\nu|$, one has
$U_{\bar{\zeta} \bar{\xi}} \subset U_\nu$ 
and $\mu(U_{\bar{\zeta} \bar{\xi}})>0$.
Since the matrix 
$A$ is irreducible and satisfies condition (I),
there exists
$\eta \in B_*(X_A)$ such that
$|\eta| >1, u \eta u \in B_*(X_A)$ and 
$U_{u \xi u} \cap U_{u \eta u} =\emptyset.
$ 
Put
$\bar{\eta} =\eta u$.
Hence
$\bar{\zeta} \bar{\eta} \in B_*(X_A)$ 
and
$U_{\bar{\zeta} \bar{\eta}} \subset U_{\nu} \subset O.
$
%Put $l = |\eta| +1$.
Define 
homeomorphisms $\psi, \varphi$ on $X_A$ by setting
for $y \in X_A$
\begin{align*}
\psi(y)
& =
{\begin{cases}
\bar{\zeta} \bar{\eta} z \in U_{\bar{\zeta}\bar{\eta}} 
& \text{ if } y = \bar{\zeta} \bar{\xi} z \in U_{\bar{\zeta} \bar{\xi}}, \\
\bar{\zeta} \bar{\xi} z \in U_{\bar{\zeta} \bar{\xi}} 
& \text{ if } y = \bar{\zeta} \bar{\eta} z \in U_{\bar{\zeta}\bar{\eta}}, \\
y
& \text{ otherwise,}
\end{cases}} \\
\varphi(y)
& =
{\begin{cases}
\bar{\zeta}\bar{\eta}\bar{\xi} z \in U_{\bar{\zeta}\bar{\eta}\bar{\xi}} 
& \text{ if } y = \bar{\zeta}\bar{\eta}\bar{\eta} z \in U_{\bar{\zeta}\bar{\eta}\bar{\eta}}, \\
\bar{\zeta}\bar{\xi}  z \in U_{\bar{\zeta}\bar{\xi}} 
& \text{ if } y = \bar{\zeta}\bar{\eta}\bar{\xi} z \in U_{\bar{\zeta}\bar{\eta}\bar{\xi}}, \\
\bar{\zeta}\bar{\eta}\bar{\eta} z \in U_{\bar{\zeta}\bar{\eta}\bar{\eta}} 
& \text{ if } y = \bar{\zeta}\bar{\xi}  z \in U_{\bar{\zeta}\bar{\xi}},\\
y
& \text{ otherwise.}
\end{cases}} 
\end{align*}
Since
$U_\zeta \subset U_\nu$,
one sees that
$$
U_{\bar{\zeta}\bar{\eta}}, 
U_{\bar{\zeta}\bar{\xi}}
U_{\bar{\zeta}\bar{\eta}\bar{\xi}}, 
U_{\bar{\zeta}\bar{\eta}\bar{\eta}} 
\subset U_\nu \subset O.
$$
Hence we have
$\psi,\varphi \in \Gamma_O$ and
$\psi^2 = \varphi^3 =\id$.
We put
$F = U_{\bar{\zeta}\bar{\xi}}\subset U_\nu $
so that 
$\mu(F) >0$ and
\begin{equation*}
\psi(F) = U_{\bar{\zeta}\bar{\eta}}, \qquad
\varphi(F) = U_{\bar{\zeta}\bar{\eta}\bar{\eta}}, \qquad
\varphi^2(F) = U_{\bar{\zeta}\bar{\eta}\bar{\xi}}.
\end{equation*} 
 Since
$$
\varphi(F) \cup \varphi^2(F) \subset \psi(F),\qquad
\varphi(F) \cap \varphi^2(F) = \emptyset,
$$
we have
$$
\mu(\varphi(F)) + \mu(\varphi^2(F)) \le \mu( \psi(F)).
$$
By hypothesis,
$\mu$ is $\Gamma_O$-invariant,
we have
$\mu(F) =0$ a contradiction.
Therefore we conclude that 
there is no $\Gamma_O$-invariant regular Borel probability measure on $O$.
\end{pf}
By using the above proof,
one may prove that the subgroup $\langle \psi, \varphi\rangle$
generated by $\psi, \varphi$ is isomorphic to
the free product ${\Bbb Z}_2 * {\Bbb Z}_3$
and hence 
$\Gamma_O$ contains the free group $F_2$ on two generators.
We will not use this fact in the discussions henceforth 
so that
we will not give its detailed proof (see \cite{MaPre2012}).

We note the following lemma.
\begin{lem} For an open set $O$ of $X_A$, 
we have  
$\gamma(O) = O$ for $\gamma \in \Gamma_O$.
\end{lem}
\begin{pf}
Take $\gamma \in \Gamma_O$.
Suppose that there is $y \in O$ such that
$\gamma(y) \in O^c$.
As $\gamma|_{O^c} = \id$,
one sees that
$\gamma(\gamma(y)) = \gamma(y)$
so that
$\gamma(y) = y \in O$ a contradiction.
 Hence we have
 $\gamma(O) \subset O$
and similarly 
 $\gamma^{-1}(O) \subset O$,
so that
$\gamma(O) = O.$
%Suppose next that  there is $x \in O$ such that  $x \in \gamma(O)^c$.
%As  $\gamma(O)^c = \gamma(O^c)$, one sees that $\gamma^{-1}(x) \in O^c$
%so that$\gamma(\gamma^{-1}(x))=\gamma^{-1}(x)$.
%This implies  $\gamma^{-1}(x) = x \in O$  a contradiction.
% Therefore we conclude that $\gamma(O) = O$.
\end{pf}
%
%\noindent
%{\bf Definition.}
The following notations also follow from \cite{GPS2}.

(1) For $O \in O(X_A)$, put $O^{\perp} =(\overline{O})^c=(O^c)^\circ \in O(X_A)$.

(2) For $F \in CL(X_A)$, put $F^{\perp} = \overline{(F^c)}=(F^\circ)^c \in CL(X_A)$.

(3) $O \in O(X_A)$ is said to be regular if $(O^{\perp})^{\perp} =O$.

(4) $F \in CL(X_A)$ is said to be regular if $(F^{\perp})^{\perp} = F$.

Denote by $RO(X_A)$ the set of regular open subsets of $X_A$.
We note that an open set $O$ of $X_A$ is regular if and only if
$O^c$ is a regular closed set of $X_A$.

\noindent
{\bf Remark.}

(i)  For an open set $O$, we have
$ (O^{\perp})^{\perp} =( ( (\overline{O})^c )^c  )^{\circ}$.
Hence $O$ is regular if and only if
$(\overline{O})^{\circ} = O $.

(ii) For a closed set $F$, we have
$ (F^{\perp})^{\perp} =\overline{( (F^{\circ})^c )^c }$.
Hence
$F$ is regular if and only if
$\overline{ F^{\circ} } = F$.

We note the following lemma.
\begin{lem}
Let $O \in O(X_A)$ be an open set.
Then 
$O^\perp$ is a regular open set.
\end{lem}
\begin{pf}
As
$({O^\perp})^\perp = (\overline{O})^\circ$,
we  have
$(({O^\perp})^\perp)^\perp =(\overline{O^\perp})^\circ \supset O^\perp$.
The inclusion relation
$(\overline{O})^\circ \supset O$
implies
$
(({O^\perp})^\perp)^\perp = ((\overline{O})^\circ)^\perp
=[ \overline{(\overline{O} )^\circ} ]^c \subset (\overline{O})^c 
= O^\perp,
$
so that we have
$(({O^\perp})^\perp)^\perp = O^\perp$.
 \end{pf}

%%%%%%%%%%%%%%%%%%%%%%%%%%%%%%%%%%%%%%%
Since
Lemma \ref{lem:3.4} is the same statement as \cite[Lemma 3.4]{GPS2},
the lemma below holds by the same proof as the proof of \cite[Lemma 3.9]{GPS2}.
\begin{lem}[cf. {\cite[Lemma 3.9]{GPS2}}] \label{lem:3.9}
For $O, O_1, O_2 \in O(X_A)$, 
we have
\begin{enumerate} 
\renewcommand{\labelenumi}{(\roman{enumi})}
\item 
 $O_1 \subset O_2 $ if and only if  $\Gamma_{O_1} \subset \Gamma_{O_2}$.
\item
 $\Gamma_O \cap \Gamma_{O^\perp} = \{\id \}.$
\item
 $(\Gamma_O)^{\perp} = \Gamma_{O^{\perp}}$ 
 and $\Gamma_O \subset  ((\Gamma_O)^{\perp})^{\perp}$.  
\item
 $O \in RO(X_A)$ 
 implies $\Gamma_O = ((\Gamma_O)^{\perp})^{\perp}$.  
\end{enumerate}
\end{lem}
\section{Strong commuting pairs}
%%%%%%%%%%%%%%%%%%%%%%%%%%%%%%%%%%%%%%%%%%%
In this section, 
we will find algebraic conditions of the pair $(\Gamma_O, \Gamma_{O^\perp})$ 
of subgroups of $\Gamma_A$ for a regular open set $O \in O(X_A)$,
and prove Proposition \ref{prop:3.13}.
   
Following \cite[Definition 3.10]{GPS2},
 we say that a pair 
$(H,K)$ of subgroups of $\Gamma_A$ is a {\it commuting pair}
if the following condition called (D1) holds: 

(D1) \, \, 
$
H^\perp = K, \qquad
K^\perp = H, \qquad
H\cap K= \{ \id \}.
$

\noindent
A commuting pair $(H,K)$ is called a {\it strong commuting pair}
if the following extra condition called $(D2)$ holds: 

(D2)
For any nontrivial normal subgroup $N$ 
of $H$ (resp. of $K$), 
then
$N^\perp = K $
(resp. 
$N^\perp = H$)
holds.

We may see similar statements in \cite{GPS2} as Lemma 3.11 and Lemma 3.12 to the following two lemmas.
The proofs of \cite[Lemma 3.11]{GPS2} and \cite[Lemma 3.12]{GPS2}
need \cite[Lemma 3.3]{GPS2}.
In our setting, we do not have a corresponding lemma to \cite[Lemma 3.3]{GPS2}.
We give complete proofs for the following two lemmas for the sake of completeness. 
\begin{lem}\label{lem:3.11}
Let $O \subset X_A$
be a nonempty open set
and
$\eta \in \Gamma_O$ with $\eta \ne \id$.
For a nonempty clopen set $U \subset O$,
there exists $\gamma\in \Gamma_O$
such that
$\gamma^{-1} \eta \gamma|_U \ne \id$.
\end{lem}
\begin{pf}
The proof below is basically similar to \cite[Lemma 3.11]{GPS2}.
Take a nonempty clopen partition
$U_1, U_2$ of $U$ so that
$U = U_1\sqcup U_2$.
Since
$\eta \ne \id$ and $\eta \in \Gamma_O$,
there exists a clopen set $Y \subset X_A$ 
such that
$Y \subset O$ and $\eta(Y) \cap Y =\emptyset$.
One may take $Y$ small enough such as
$U_1 \backslash \eta(Y) \ne \emptyset$
and
$U_2 \backslash Y \ne \emptyset.$
By Lemma \ref{lem:2.1}, 
there exist 
$\alpha \in \Gamma_A$ and
a nonempty clopen set $U'_1 \subset U_1 \backslash \eta(Y)$
such that
\begin{equation*}
\alpha(U'_1) \subset Y,\qquad
\alpha^2 = \id,\qquad
\alpha|_{(U'_1\cup\alpha(U'_1))^c} = \id.
\end{equation*}
For $U_2\backslash Y$ and 
$\eta(\alpha(U'_1))(\subset \eta(Y))$,
Lemma \ref{lem:2.1} assures that there exist
$\beta \in \Gamma_A$ and
a nonempty clopen set $U'_2 \subset U_2 \backslash Y$
such that
\begin{equation*}
\beta(U'_2) \subset \eta(\alpha(U'_1)),\qquad
\beta^2 = \id,\qquad
\beta|_{(U'_2\cup\beta(U'_2))^c} = \id.
\end{equation*}
As 
$\alpha(U'_1) \cap U'_2 \subset Y\cap U'_2 = \emptyset$,
$\beta(U'_2) \cap U'_1 \subset \eta(Y)\cap U'_1 = \emptyset$,
one notes that
$$
[U'_1\cup\alpha(U'_1)] \cap [U'_2\cup\beta(U'_2)] = \emptyset.
$$
Define
$\gamma \in \Gamma_A$ by setting
\begin{equation*}
\gamma
= 
\begin{cases}
\alpha & \text{ on } U'_1\cup\alpha(U'_1),\\
\beta & \text{ on } U'_2\cup\beta(U'_2),\\
\id   & \text{ elsewhere}.
\end{cases}
\end{equation*}
As
$U'_1 \subset U_1 \subset O$
and
$\alpha(U'_1) \subset Y \subset O$,
and also
$U'_2 \subset U_2 \subset O$
and
$\beta(U'_2) \subset 
 \eta(Y)\subset  O$,
both
$U'_1\cup\alpha(U'_1) $ 
and
$U'_2\cup\beta(U'_2) $
are subsets of $O$.
Hence 
$\gamma $ belongs to $\Gamma_O$.
Since
$\beta(U'_2) \subset 
\eta(\alpha(U'_1)),
$
we have
$$
\alpha \eta^{-1}\beta(U'_2) \subset 
\alpha^2(U'_1) = U'_1
$$
so that
$$
\gamma^{-1}\eta \gamma (\alpha \eta^{-1}\beta(U'_2)) 
=
\gamma^{-1}\eta \alpha\alpha \eta^{-1}\beta(U'_2) 
=
\gamma^{-1}\beta(U'_2) 
= U'_2.
$$
Therefore we have
$\gamma^{-1} \eta \gamma|_U \ne \id$.
\end{pf} 
\begin{lem}\label{lem:3.12}
Let $O \subset X_A$
be a nonempty open set
and
$\eta \in \Gamma_O$ with $\eta \ne \id$.
For a nonempty clopen set $U \subset O$
and
$\gamma \in \Gamma_A$ 
such that
$\gamma(U) \subset O$ and
$U \cap \gamma(U) =\emptyset$,
there exist
a nonempty clopen set
$U_1 \subset U$ and
 $\psi\in \Gamma_O$
such that
\begin{equation*}
\gamma(\psi^{-1} \eta \psi)(U_1) \cap
(\psi^{-1} \eta \psi)(\gamma(U_1))
=\emptyset.
\end{equation*}
\end{lem}
\begin{pf}
As
$U \cup \gamma(U) \subset O$,
by taking a smaller clopen set $U'\subset U$
such as
$U' \cap \gamma(U') \ne O$,
one may assume that
$U \cup \gamma(U)$ is a proper subset of $O$
and take  a nonempty clopen subset 
$\tilde{U} \subset O\backslash{(U \cup\gamma(U))}$.
For the clopen set 
$\tilde{U}$ and $\eta \in \Gamma_O$ with $\eta \ne \id$,
by applying the proof of the preceding lemma,
we have a nontrivial clopen partition 
$\tilde{U}_1 \sqcup\tilde{U}_2 =\tilde{U}$,
a nonempty clopen subset $\tilde{U}'_2$ of $\tilde{U}_2$
and
$\tilde{\gamma} \in \Gamma_A$ 
such that
$$
\tilde{\gamma}^{-1}\eta \tilde{\gamma} 
(\alpha \eta^{-1}\beta(\tilde{U}'_2))
\subset
\tilde{U}'_2 \subset \tilde{U}_2.
$$
Put
$Y = \alpha \eta^{-1}\beta(\tilde{U}'_2) \subset \tilde{U}_1$ 
so that
$
\tilde{\gamma}^{-1}\eta \tilde{\gamma}(Y)
\subset
\tilde{U}_2
$
and
$
\tilde{\gamma}^{-1}\eta \tilde{\gamma}(Y)
\cap
Y = \emptyset.
$
By putting  
$
\tilde{\gamma}^{-1}\eta \tilde{\gamma}
$
as
$\eta $,
one has
%a nonempty clopen set $Y \subset X_A$ such that
$$
Y \cap \eta(Y) = \emptyset, \qquad 
Y \cup \eta(Y) \subset O\backslash{(U\cup\gamma(U))}.
$$
Let
$
U = U^{(1)} \sqcup U^{(2)} \sqcup U^{(3)}
$
be a nontrivial clopen partition of $U$.
By Lemma \ref{lem:2.1},
there exist
$\hat{\alpha}\in \Gamma_A
$ 
with
a clopen subset $\widehat{U}_1 \subset U^{(1)}$
such that
$\hat{\alpha}(\widehat{U}_1) \subset U^{(2)}$
and
$\hat{\beta}\in \Gamma_A
$
with a clopen subset
$\widehat{U}_2 \subset \hat{\alpha}(\widehat{U}_1)$
such that
$\hat{\beta}(\widehat{U}_2) \subset U^{(3)}$.
Put
\begin{equation*}
U' = \hat{\alpha}^{-1}(\widehat{U}_2) \subset \widehat{U}_1 \subset  U^{(1)},
\qquad
U^{\prime\prime} = \widehat{U}_2 \subset U^{(2)},
\qquad
U^{\prime\prime\prime} 
=\hat{\beta}(\widehat{U}_2) \subset U^{(3)}.
\end{equation*}
Hence we have a disjoint $\Gamma_A$-equivalent clopen partition:
$$
U^{\prime}\sqcup U^{\prime\prime} \sqcup
U^{\prime\prime\prime} \subset U.
$$
Let
$Y = Y_1 \sqcup Y_2$
be a nontrivial clopen partition of $Y$.
By Lemma \ref{lem:2.1}, there exist 
$U_1 \subset U'$
and
$\alpha,\beta \in \Gamma_A$
such that
\begin{align}
&\alpha^2 =  \beta^2 = \id, \qquad
\alpha(U_1) \subset Y_1, \qquad
\alpha(\gamma(U_1)) \subset Y_2, \\
&\beta(\eta\alpha(U_1)))  \subset U^{\prime\prime}, \qquad
\beta(\eta\alpha\gamma(U_1)))  \subset \gamma(U^{\prime\prime\prime})  \label{eqn:beta1} \\
\intertext{and}
&\alpha|_{[U_1 \cup \alpha(U_1)\cup \gamma(U_1) \cup \alpha\gamma(U_1)]^c}
=\id, \\
& \beta|_{[\eta\alpha(U_1)\cup \beta\eta\alpha(U_1)
\cup \eta\alpha\gamma(U_1) \cup \beta\eta\alpha\gamma(U_1)]^c}
=\id. \label{eqn:beta2}
\end{align}
In fact, by applying Lemma \ref{lem:2.1} for $U'$ and $Y_1$,
we have  $V_1 \subset U'$ and $\alpha_1 \in \Gamma_A$ such that
\begin{equation*}
\alpha_1(V_1) \subset Y_1, \quad
\alpha_1^2 = \id, \quad
\alpha_1|_{(V_1 \cup \alpha_1(V_1))^c} =\id.
\end{equation*}
By applying Lemma \ref{lem:2.1} for $\gamma(V_1)$ and $Y_2$,
we have  $V'_1 \subset \gamma(V_1)$ and $\alpha_2 \in \Gamma_A$ such that
\begin{equation*}
\alpha_2(V'_1) \subset Y_2, \quad
\alpha_2^2 = \id, \quad
\alpha_2|_{(V'_1 \cup \alpha_2(V'_1))^c} =\id.
\end{equation*}
Put
$U_1 =\gamma^{-1}(V'_1)$
so that
$U_1 \subset V_1$
and
\begin{equation*}
\alpha_1(U_1) \subset Y_1, \quad
\alpha_1|_{(U_1 \cup \alpha_1(U_1))^c} =\id, \quad
\alpha_2\gamma(U_1) \subset Y_2, \quad
\alpha_2|_{ ( \gamma(U_1) \cup \alpha_2\gamma(U_1) )^c} =\id.
\end{equation*}
Since
%$\gamma(U_1) \cup \alpha_2(\gamma(U_1)) \subset \gamma(U) \cup Y_2$ and
%$U_1 \cup \alpha_1(U_1) \subset U \cup Y_1$.
%As $\gamma(U) \subset U^c$ and $Y_2 \subset Y_1^c$,
%we have $\gamma(U_1) \cup \alpha_2(\gamma(U_1)) \subset 
$
[U_1 \cup \alpha_1(U_1)\cup \gamma(U_1) \cup \alpha_2\gamma(U_1)]^c
\subset
(U_1 \cup \gamma(U_1) )^c 
$
and
$U_1 \cap \gamma(U_1) =\emptyset$,
by putting
\begin{equation*}
\alpha =
\begin{cases}
\alpha_1 & \text{ on } U_1, \\
\alpha_2 & \text{ on } \gamma(U_1), \\
\id & \text{ elsewhere,}
\end{cases}
\end{equation*}
we have
\begin{equation*}
\alpha^2 =  \id, \quad
\alpha(U_1) \subset Y_1, \quad
\alpha(\gamma(U_1)) \subset Y_2, \quad
\alpha|_{[U_1 \cup \alpha(U_1)\cup \gamma(U_1) \cup \alpha\gamma(U_1)]^c}
=\id.
\end{equation*}
By a similar manner to the above discussion,
 we have $\beta \in \Gamma_A$ with $\beta^2 = \id$
satisfying
\eqref{eqn:beta1} and \eqref{eqn:beta2}.
We define
$\psi \in \Gamma_A$
by setting
\begin{equation*}
\psi =
\begin{cases}
\alpha & \text{ on } U_1 \cup \alpha(U_1)\cup \gamma(U_1) \cup \alpha\gamma(U_1), \\
\beta & \text{ on } 
\eta\alpha(U_1) \cup \beta\eta\alpha(U_1)
\cup \eta\alpha\gamma(U_1) \cup \beta\eta\alpha\gamma(U_1),
\\
\id & \text{ elsewhere.}
\end{cases}
\end{equation*}
We then have
$$
\psi^{-1}\eta\psi(U_1) \subset U^{\prime\prime} \subset U,
\qquad
\psi^{-1}\eta\psi(\gamma(U_1)) \subset \gamma(U^{\prime\prime\prime})
 \subset \gamma(U)
$$
so that
\begin{equation*}
\gamma(\psi^{-1}\eta\psi(U_1)) \cap
(\psi^{-1}\eta\psi)(\gamma(U_1)) \
\subset \
\gamma(U^{\prime\prime}) \cap \gamma(U^{\prime\prime\prime})
=
\emptyset.
\end{equation*}
\end{pf}
By  
using Lemma \ref{lem:3.11}
and Lemma \ref{lem:3.12}
with Lemma 2.6, 
one may show the following proposition
 in a similar manner to
the proof of \cite[Proposition 3.13]{GPS2}.
\begin{prop} \label{prop:3.13}
If $O $ is a regular open set of $X_A$,
then the pair
$(\Gamma_O, \Gamma_{O^\perp})$ 
is a strong commuting pair.
\end{prop}
\section{The clopen condition}
%%%%%%%%%%%%%%%%%%%%%%%%%%%%%%%%%%%
In this section, 
we will find algebraic conditions of the strong commuting pair 
$(\Gamma_O, \Gamma_{O^\perp})$ 
of subgroups of $\Gamma_A$ for a clopen set $O \in CO(X_A)$,
and prove Proposition \ref{prop:4.12}.
The proposition comes from Lemma \ref{lem:3.21},
which is based on Lemma \ref{lem:3.18} and Lemma \ref{lem:3.19}.
There are similar statements to Lemma \ref{lem:3.18} and Lemma \ref{lem:3.19}
in \cite{GPS2} as \cite[Lemma 3.18]{GPS2} and \cite[Lemma 3.19]{GPS2}.
In their proofs of the latter lemmas, \cite[Lemma 3.3]{GPS2} has been used.
In our setting, we do not have a version of \cite[Lemma 3.3]{GPS2}, 
so that we must modify the given proofs in \cite{GPS2}
of \cite[Lemma 3.18]{GPS2} and \cite[Lemma 3.19]{GPS2}. 
In particular, we must provide several lemmas to prove Lemma \ref{lem:3.18}.

\begin{lem}\label{lem:4.1}
For a word $\nu=(\nu_1,\dots,\nu_n) \in B_n(X_A)$ with $n >1$
and a nonempty open set $V \subset X_A$ such that $U_\nu$ does not 
contain $V$,
there exists $\alpha \in \Gamma_A$ 
such that
\begin{equation*}
\alpha(U_\nu) \subset V, 
\qquad
\alpha^2 = \id,
\qquad
\alpha|_{(U_\nu \cup \alpha(U_\nu))^c} =\id.
\end{equation*}
\end{lem}
\begin{pf}
Take $\mu = \mu_1\cdots\mu_k\in B_k(X_A)$
such that
$k >n$ and 
$(\mu_1,\dots,\mu_n) \ne (\nu_1,\dots,\nu_n)$.
As $A$ is irreducible,
there exists a word
$(\xi_1,\dots,\xi_l )\in B_l(X_A)$
such that
$\mu \xi \nu \in B_*(X_A)$.
Hence we have
$$
U_{\mu \xi \nu} \subset U_\mu \subset V,
\qquad
U_{\mu \xi \nu} \cap U_\nu =\emptyset.
$$
Define 
$\alpha \in \Gamma_A$ 
by setting for $x \in X_A$
\begin{align*}
& \alpha(x) \\ 
=&
{\begin{cases}
& (\mu_1, \dots,\mu_k, \xi_1,\dots,\xi_l,\nu_1,\dots,\nu_n,
x_{n+1},x_{n+2},\dots) \in U_{\mu \xi\nu}, \\ 
& \text{ if } x=(\nu_1,\dots,\nu_n,x_{n+1},x_{n+2},\dots)\in U_\nu, \\
& (\nu_1, \dots,\nu_n,
x_{n+1},x_{n+2},\dots) \in U_{\nu}, \\
& \text{ if } x=
(\mu_1,\dots,\mu_k, \xi_1,\dots,\xi_l,\nu_1,\dots,\nu_n, x_{n+1},x_{n+2},\dots)\in U_{\mu \xi \nu}, \\
& x  \quad \text{ otherwise.}
\end{cases}} 
\end{align*}
Then $\alpha$ defines an element of $\Gamma_A$
which has the desired properties.
\end{pf}
\begin{lem}\label{lem:4.2}
For $U, V \in CO(X_A)$ with
$V \backslash U \ne \emptyset$,
there exist a clopen partition
$U_1,\dots,U_n \subset U$ of $U$
such that
$\cup_{i=1}^n U_i =U$ and 
$U_i \cap U_j =\emptyset$ for $i\ne j$
and homeomorphisms
$\alpha_1,\dots,\alpha_n \in \Gamma_A$
such that
\begin{equation*}
\alpha_i(U_i) \subset V\backslash U, \qquad
\alpha_i(U_i) \cap \alpha_j(U_j) =\emptyset,  \qquad
\alpha_i^2 = \id, \qquad
{\alpha_i|}_{(U_i \cup \alpha(U_i))^c} =\id
\end{equation*}
for $i,j=1,\dots, n$ with $i\ne j$.
 \end{lem}
\begin{pf}
Since $U$ is clopen,
there exist words 
$\nu(1), \dots,\nu(n) \in B_*(X_A)$
such that $U$ is a disjoint union of cylinder sets
$U = U_{\nu(1)}\sqcup \cdots\sqcup U_{\nu(n)}$.
Put $\widetilde{V} = V \backslash U$ a nonempty clopen set.
Take $V_1 \subset \widetilde{V}$ a clopen subset of $\widetilde{V}$
such that
$V_1 \ne \widetilde{V}$.
For $U_{\nu(1)}$ and $V_1$, Lemma \ref{lem:4.1} 
ensures that
there exists 
$\alpha_1 \in \Gamma_A$ such that
\begin{equation*}
\alpha_1(U_{\nu(1)}) \subset V_1, \qquad
\alpha_1^2 = \id, \qquad
{\alpha_1}|_{(U_{\nu(1)} \cup \alpha_1(U_{\nu(1)}))^c} =\id.
\end{equation*}
By applying Lemma \ref{lem:4.1},  recursively,
 we have  clopen sets
 $V_1,\dots,V_n \subset \widetilde{V}$ with
 $V_i \cap V_j =\emptyset $ for
 $i\ne j$ and
$\alpha_i \in \Gamma_A, i=1,\dots,n$ 
such that
\begin{equation*}
\alpha_i(U_{\nu(i)}) \subset V_i, \qquad
\alpha_i^2 = \id, \qquad
{\alpha_i}|_{(U_{\nu(i)} \cup \alpha_i(U_{\nu(i)}))^c} =\id
\end{equation*}
for
$i=1,\dots,n$.
By putting $U_i= U_{\nu(i)}$,
the proof ends.
\end{pf}
\begin{lem}\label{lem:4.3}
Let $U, W \in CO(X_A)$ be nonempty clopen sets such that
$U \cap W = \emptyset$.
Then there exists
$\alpha \in \Gamma_A$
such that
\begin{equation*}
\alpha(U) \subset W, \qquad
\alpha^2 = \id, \qquad
\alpha|_{(U \cup \alpha(U))^c} =\id.
\end{equation*}\end{lem}
\begin{pf}
By the preceding lemma, we have clopen partitions:
$U_1\sqcup \cdots \sqcup U_n = U,$
$W_1\sqcup \cdots \sqcup W_n = W$
and
homeomorphisms
$\alpha_i \in \Gamma_A, i=1,\dots,n$
such that  
\begin{equation*}
\alpha_i(U_i) \subset W_i, \qquad
\alpha_i^2 = \id, \qquad
{\alpha_i}|_{(U_{i} \cup \alpha_i(U_{i}))^c} =\id
\end{equation*}
for
$i=1,\dots,n$.
The homeomorphisms 
$\alpha_i, i=1,\dots,n$ commute with each other.
The homeomorphism
$\alpha = \alpha_1\circ \cdots \circ \alpha_n \in \Gamma_A$
has the desired properties.
\end{pf}

\begin{lem}\label{lem:4.4}
Let $O \in CO(X_A)$ be a clopen set.
Let
$U \subset O, \ V \subset O^c$ be $\Gamma_A$-equivalent nonempty clopen sets
and
$W\subset O, \ W'\subset O^c$ be nonempty clopen sets
such that
$U \cap W =\emptyset, V \cap W' =\emptyset$.
Then
there exist a clopen partition
$U_1,\dots,U_n \subset U$ of $U$ 
and a clopen partition
$V_1,\dots,V_n \subset V$ of $V$
such that
$U = \sqcup_{i=1}^n U_i, V = \sqcup_{i=1}^n V_i,$ 
$U_i \underset{\Gamma_A}{\sim} V_i$ for $i=1,\dots,n$
and homeomorphisms
$\alpha_i \in \Gamma_O, \ \beta_i \in \Gamma_{O^\perp}$
for $i=1,\dots,n$ 
such that
\begin{align}
\alpha_i(U_i) & \subset W, \quad
\alpha_i(U_i) \cap \alpha_j(U_j) =\emptyset, \quad
\alpha_i^2 = \id, \quad
{\alpha_i|}_{(U_i \cup \alpha(U_i))^c} =\id, \label{eqn:alphau} \\ 
\beta_i(V_i) & \subset W', \quad
\beta_i(V_i) \cap \beta_j(V_j) =\emptyset, \quad
\beta_i^2 = \id, \quad
{\beta_i|}_{(V_i \cup \beta(V_i))^c} =\id \label{eqn:betav}
\end{align}
for $i,j=1,\dots, n$ with $i\ne j$.
 \end{lem}
\begin{pf}
Since $U \underset{\Gamma_A}{\sim} V$,
there exists
$\gamma \in \Gamma_A$ such that 
$\gamma(U) = V$.
As $U \cap V =\emptyset$,
by Lemma \ref{lem:2.2}, 
one may assume that 
$\gamma^2 = \id$ and
$\gamma|_{(U\cup V)^c} =\id$.
Since $\gamma$ is a cylinder map,
there exist words
$\nu(1),\dots,\nu(n),\mu(1),\dots,\mu(n)\in B_*(X_A)$
such that
\begin{gather*}
U=\cup_{i=1}^n U_{\nu(i)},\qquad
V=\cup_{i=1}^n U_{\mu(i)},\qquad 
\gamma(U_{\nu(i)}) = U_{\mu(i)}, \quad i=1,\dots,n, \\
U_{\nu(i)}\cap U_{\nu(j)} = \emptyset, \qquad
U_{\mu(i)}\cap U_{\mu(j)}=\emptyset, \quad
i\ne j.
\end{gather*}
Hence
$U_{\nu(i)} \underset{\Gamma_A}{\sim} V_{\mu(i)}, \ i=1,\dots,n.$
By the proof of  Lemma \ref{lem:4.2},
one may find
$\alpha_i \in \Gamma_A$ 
such that
\begin{align*}
\alpha_i(U_{\nu(i)}) & \subset W, \qquad
\alpha_i(U_{\nu(i)}) \cap \alpha_j(U_{\nu(j)}) =\emptyset,  \\
\alpha_i^2 & = \id, \qquad
{\alpha_i |}_{(U_{\nu(i)} \cup \alpha(U_{\nu(i)}))^c}   =\id
\end{align*}
for $i,j=1,\dots,n$ with $i\ne j$.
As $U_{\nu(i)} \cup \alpha_i(U_{\nu(i)}) \subset O$,
one sees that $\alpha_i \in \Gamma_O$.
One may similarly find $\beta_i \in \Gamma_{O^\perp}$
having the desired properties
by putting 
$U_i = U_{\nu(i)}, \ V_i = U_{\mu(i)}$
for
$i=1,\dots,n$.
\end{pf}
For  subsets $H_1, H_2, H_3$ of $\Gamma_A$,
let us denote by $\langle H_1, H_2 \rangle$ and $\langle H_1, H_2, H_3 \rangle$ 
the subgroup of $\Gamma_A$
generated by  elements of $H_1, H_2$ and 
that of $\Gamma_A$
generated by  elements of $H_1, H_2, H_3$
respectively.

One may see a similar statement to 
the following two lemmas in \cite[Lemma 3.18]{GPS2}.
The proof of \cite[Lemma 3.18]{GPS2} needs \cite[Lemma 3.3]{GPS2} 
for which we do not have a corresponding lemma in our setting.
The following proofs are different from the proof of \cite[Lemma 3.18]{GPS2}.
\begin{lem}\label{lem:4.5}
Let $O \in CO(X_A)$ 
and $\eta \in \Gamma_A$ 
satisfy
$
\eta(O)\cap O^c \ne \emptyset,
\eta(O^c)\cap O^c \ne \emptyset.
$
Let
$U \subset O, \ V \subset O^c$ be $\Gamma_A$-equivalent nonempty clopen sets
such that
$
O\cap \eta^{-1}(O^c) \cap U^c   \ne \emptyset,
O^c\cap \eta^{-1}(O^c) \cap V^c \ne \emptyset.
$
Then
there exists 
$\chi \in \langle\Gamma_O,\Gamma_{O^\perp},\eta\rangle$
such that
\begin{equation*}
\chi(U) = V,\qquad
\chi(V) =U, \qquad
\chi|_{(U \cup V)^c} =\id.
\end{equation*}
 \end{lem}
\begin{pf}
Put the nonempty clopen sets
$$
W = O\cap \eta^{-1}(O^c) \cap U^c, \qquad
W' =O^c\cap \eta^{-1}(O^c) \cap V^c.
$$
By the preceding lemma,
there exist disjoint clopen partitions of $U$ and of $V$:
$$
U = U_1\sqcup \cdots \sqcup U_n, 
\qquad 
V = V_1\sqcup \cdots \sqcup  V_n
$$
such that
$U_i \underset{\Gamma_A}{\sim} V_i$ 
and
$\alpha_i \in \Gamma_O, \ 
\beta_i \in \Gamma_{O^\perp}$ 
satisfying \eqref{eqn:alphau} and \eqref{eqn:betav}
for
$i=1,\dots,n$.
Put
$U^i = \alpha_i(U_i) \subset W, \
V^i = \beta_i(V_i) \subset W'.
$
One then sees 
$$
\eta(U^i), \eta(V^i)
\subset O^c,
\qquad
\eta(U^i) \underset{\Gamma_A}{\sim} \eta(V^i),
\qquad
\eta(U^i) \cap \eta(V^i) =\emptyset, \qquad i=1,\dots,n.
$$
By Lemma \ref{lem:2.2},
there exist
$\gamma^i \in \Gamma_{O^\perp}, i=1,\dots,n$
such that
$$
\gamma^i(\eta(U^i)) = \eta(V^i), \qquad
{(\gamma^i)}^2 = \id, \qquad
\gamma^i|_{(\eta(U^i)\cup \eta(V^i))^c}= \id.
$$
Then
$
\chi^i = \eta^{-1}\gamma^i\eta \in \langle\Gamma_{O^\perp},\eta\rangle 
$
satisfies
$$
\chi^i(U^i) = V^i, \qquad \chi^i(V^i) = U^i,
\qquad i=1,\dots,n.
$$
For $x \in (U^i \cup V^i)^c$,
we have
%$\eta(x) \in (\eta(U^i) \cup \eta(V^i))^c$ so that
$\gamma^i(\eta(x) ) = \eta(x)$ 
and hence
$\chi^i|_{(U^i \cup V^i)^c} =\id$.
We put
$$
\chi_i = \alpha_i \beta_i \chi^i \alpha_i \beta_i  
\in 
\Gamma_O \Gamma_{O^\perp} \langle\Gamma_{O^\perp},\eta\rangle\Gamma_O \Gamma_{O^\perp},
\qquad
i=1,\dots,n.
$$
It then follows that
\begin{align*}
\chi_i(U_i)
& = \alpha_i \beta_i \chi^i \alpha_i (U_i) 
 = \alpha_i \beta_i \chi^i(U^i)
 = \alpha_i \beta_i (V^i)
 = \alpha_i (V_i) = V_i, \\
\chi_i(V_i)
& = \alpha_i \beta_i \chi^i \alpha_i (V_i) 
 = \alpha_i \beta_i \chi^i(V^i)
 = \alpha_i \beta_i (U^i)
 = \alpha_i (U^i) = U_i,  \\
(U_i \cup V_i)^c 
& =(\alpha_i(U^i) \cup \beta_i(V^i))^c
=(\alpha_i\beta_i(U^i) \cup \alpha_i\beta_i(V^i))^c
=\alpha_i\beta_i((U^i \cup V^i)^c).
\end{align*}
As
$
\alpha_i \beta_i \chi^i \alpha_i \beta_i 
|_{\alpha_i\beta_i((U^i \cup V^i)^c)} =\id,
$
one sees 
$\chi_i|_{(U_i \cup V_i)^c} = \id.
$
We set
$$
\chi =\chi_1\chi_2\cdots \chi_n \in \langle\Gamma_O,\Gamma_{O^\perp},\eta\rangle.
$$
Since
$
\chi_i \chi_j =  \chi_j \chi_i
$ 
for $i,j=1,\dots,n$,
we have
\begin{equation*}
\chi(U) = V,\qquad
\chi(V) =U, \qquad
\chi|_{(U \cup V)^c} =\id.
\end{equation*}
\end{pf}

\begin{lem}\label{lem:3.18}
Let $O \in CO(X_A)$ 
and $\eta \in \Gamma_A$ 
satisfy
$
\eta(O)\cap O^c \ne \emptyset,
\eta(O^c)\cap O^c \ne \emptyset.
$
Let
$U \subset O, \ V \subset O^c$ be $\Gamma_A$-equivalent nonempty clopen sets.
Then
there exists 
$\chi \in \langle\Gamma_O,\Gamma_{O^\perp}, \eta\rangle$
such that
\begin{equation*}
\chi(U) = V,\qquad
\chi(V) =U, \qquad
\chi|_{(U \cup V)^c} =\id.
\end{equation*}
 \end{lem}
\begin{pf}
We first assume that
$ U \ne  O$ and $V \ne O^c$.
If both of the  conditions
$
O\cap \eta^{-1}(O^c) \cap U^c   \ne \emptyset
$ 
and
$ 
O^c\cap \eta^{-1}(O^c) \cap V^c \ne \emptyset
$
hold,
then the assertion follows from the previous lemma.
We next assume that
\begin{equation*}
O\cap \eta^{-1}(O^c) \cap U^c   = \emptyset, \qquad
O^c\cap \eta^{-1}(O^c) \cap V^c = \emptyset
\end{equation*}
and hence
\begin{equation*}
O\cap \eta^{-1}(O^c) \subset U, \qquad
O^c\cap \eta^{-1}(O^c) \subset  V.
\end{equation*}
By Lemma \ref{lem:4.3},
there exist
 clopen sets
$U' \subset O$ and 
$V' \subset O^c$,
and homeomorphisms
$\alpha \in \Gamma_O$ and
$\beta \in \Gamma_{O^\perp}$
such that
\begin{align}
U' \cap U & = \emptyset, \quad
\alpha(U) =U', \quad
\alpha^2 =\id,\quad
\alpha|_{(U \cup U')^c} =\id, \label{eqn:u'u}\\
\intertext{and}
V' \cap V & = \emptyset, \quad
\beta(V) =V', \quad
\beta^2 =\id,\quad
\beta|_{(V \cup V')^c} =\id. \label{eqn:v'v}
\end{align}
Since
$U'$ and $V'$ are $\Gamma_A$-equivalent,
they
satisfy the assumption of the preceding lemma,
so that
there exists
$\widetilde{\chi} \in 
\langle \Gamma_O, \Gamma_{O^\perp}, \eta\rangle
$
such that
\begin{equation*}
\widetilde{\chi}(U') = V',\qquad
\widetilde{\chi}(V') = U', \qquad
\widetilde{\chi}|_{(U' \cup V')^c} =\id.
\end{equation*}
We then see that
the homeomorphism
$\chi = \alpha \circ \beta \circ \widetilde{\chi}\circ \alpha \circ \beta$
belongs to
$\langle \Gamma_O, \Gamma_{O^\perp}, \eta\rangle.
$
As
$\beta(U) =U, \beta(V') = V, \alpha(V) =V$,
it then follows that
\begin{equation*}
\chi(U) 
=\alpha \circ \beta \circ \widetilde{\chi}\circ \alpha(U)
=\alpha \circ \beta \circ \widetilde{\chi}(U')
=\alpha \circ \beta (V')
=\alpha(V)
=V,
\end{equation*}
and similarly
$\chi(V) = U$.
We note that
$\alpha$ commutes with $\beta$.
For $x \in (U\cup V)^c$,
we have
$\alpha(x) \in (U')^c, \beta(x) \in (V')^c$
and
$$
\alpha\beta(x) \in \alpha((V')^c) = (\alpha(V'))^c =(V')^c,\quad
\beta\alpha(x) \in \beta((U')^c) = (\beta(U'))^c = (U')^c
$$
so that
$\alpha\beta(x) \in (U' \cup V')^c$.
As
$\widetilde{\chi}|_{(U' \cup V')^c} =\id$,
one then sees 
\begin{equation*}
\chi(x) = \alpha\beta\widetilde{\chi}\alpha\beta(x)
= \alpha\beta\alpha\beta(x)
=x.
\end{equation*}
This shows that
$\chi|_{(U \cup V)^c} =\id$.

If 
$
O\cap \eta^{-1}(O^c) \cap U^c   = \emptyset
$ 
and
$ 
O^c\cap \eta^{-1}(O^c) \cap V^c \ne \emptyset,
$
then we may take $U' \subset O$ instead of $U$
such that \eqref{eqn:u'u},
 and apply the preceding lemma for
 $U' \subset O$  and $V \subset O$.
We then have 
$\chi' \in \langle \Gamma_O, \Gamma_{O^\perp}, \eta \rangle$
such that
\begin{equation*}
{\chi'}(U') = V,\qquad
{\chi'}(V) = U', \qquad
{\chi'}|_{(U' \cup V)^c} =\id.
\end{equation*}
By putting 
$\chi = \alpha \circ \chi' \circ \alpha$,
we have a desired homeomorphism.

If 
$
O\cap \eta^{-1}(O^c) \cap U^c   \ne \emptyset
$ 
and
$ 
O^c\cap \eta^{-1}(O^c) \cap V^c = \emptyset,
$
we symmetrically have a desired homeomorphism.

We will finally consider general clopen sets
$U \subset O$ and
$V \subset O^c$.
The conditions
$U \ne O, V \ne O^c$ are not necessarily assumed.
Suppose that $U$ and $V$ are $\Gamma_A$-equivalent
so that there exists $\gamma\in \Gamma_A$ such that 
$\gamma(U) = V$.
Take nonempty clopen sets
$U^1, U^2$ such that
$U = U^1 \cup U^2$ and $U^1 \cap U^2 =\emptyset$.
Put
$V^1 =\gamma(U^1), V^2 =\gamma(U^2)$
Hence
$U_i \underset{\Gamma_A}{\sim} V_i$.
Since
$U_i \ne O$, $V_i \ne O^c$,
by the above discussions,
one may find
$\chi_i \in \langle \Gamma_O,\Gamma_{O^\perp},\eta \rangle$
such that
\begin{equation*}
\chi_i(U_i) = V_i,\qquad
\chi_i(V_i) = U_i, \qquad
\chi_i|_{(U_i \cup V_i)^c} =\id,
\qquad i=1,2.
\end{equation*}
As $\chi_1\chi_2 = \chi_2\chi_1$,
by setting
$\chi = \chi_1\circ \chi_2 
\in  \langle \Gamma_O,\Gamma_{O^\perp},\eta \rangle$,
we have
\begin{equation*}
\chi(U) = V,\qquad
\chi(V) = U, \qquad
\chi|_{(U \cup V)^c} =\id.
\end{equation*}
\end{pf}
\begin{lem}
Let $O \in CL(X_A), \gamma \in \Gamma_A$ satisfy
$\gamma(O) = O$,
then there exist 
$\gamma_1 \in \Gamma_O, \gamma_2 \in \Gamma_{O^\perp}$
such that
$\gamma = \gamma_1 \gamma_2$
so that
$\gamma \in \Gamma_O \Gamma_{O^\perp}.
$
\end{lem}
\begin{pf}
Assume that
$\gamma(O) = O$.
We set for $ x \in \Gamma_A$
\begin{equation*}
\gamma_1(x)  =
{\begin{cases}
\gamma(x) & \text{ if } x \in O, \\
x         & \text{ if } x \in O^c,
\end{cases}},  \qquad
\gamma_2(x)  =
{\begin{cases}
x          & \text{ if } x \in O, \\
\gamma(x)  & \text{ if } x \in O^c.
\end{cases}}
\end{equation*}
The homeomorphisms
$\gamma_1, \gamma_2$ satisfy
$\gamma_1 \in \Gamma_O, \gamma_2 \in \Gamma_{O^\perp}$
and
$\gamma = \gamma_1 \gamma_2$.
\end{pf}
The same statement as the following lemma is seen in 
\cite[Lemma 3.19]{GPS2}.
The proof of \cite[Lemma 3.19]{GPS2} 
however 
does not well work in our setting.
We give a different proof as in the following way.   
\begin{lem}\label{lem:3.19}
Let $O \in CL(X_A)$ and $\eta \in \Gamma_A$
satisfy
$\eta(O) \cap O^c \ne \emptyset,
\eta(O^c) \cap O^c \ne \emptyset.
$
Then the subgroup
$\langle \Gamma_O, \Gamma_{O^\perp}, \eta\rangle
$ 
coincides with
$\Gamma_A$.
\end{lem}
\begin{pf}
For an arbitrary fixed homeomorphism $\psi \in \Gamma_A$,
we will show that
$
\psi \in \langle \Gamma_O, \Gamma_{O^\perp}, \eta\rangle.
$ 
By Lemma \ref{lem:4.3},
there exist $\alpha_1 \in \Gamma_A$ such that 
$\alpha_1(O) \subset O^c, 
\alpha_1^2 = \id, \alpha_1 |_{(O\cup \alpha(O))^c}=\id.
$
Put
$V_1 =\alpha_1(O)$.
By Lemma \ref{lem:3.18},
there exists
$\chi_1 \in
\langle \Gamma_O, \Gamma_{O^\perp}, \eta\rangle
$
such that
$\chi_1(O) = V_1, \chi_1(V_1) = O, \chi_1|_{(O \cup V_1)^c} = \id.$
We have two cases.

Case 1: $O \cap \psi^{-1}(O^c) \ne \emptyset.$

Since
$O \cap \psi^{-1}(O^c) \subset O$
and
$V_1 \subset O^c$,
by Lemma \ref{lem:4.3}, 
there exists
$\alpha_2 \in \Gamma_A$ 
such that
$\alpha_2(V_1) \subset O\cap \psi^{-1}(O^c)$.
Put
$V_2 =\alpha_2(V_1)$.
By Lemma \ref{lem:3.18}, 
one may take 
$\chi_2\in 
\langle \Gamma_O, \Gamma_{O^\perp}, \eta\rangle
$
such that
$\chi_2(V_1) = V_2$.
Put
$W_2 = \psi(V_2) \subset O^c$
so that
$\psi\chi_2\chi_1(O) = W_2 \subset O^c$.
Since
$O \underset{\Gamma_A}{\sim}W_2$,
by Lemma \ref{lem:3.18}, 
we have 
$\chi_3\in 
\langle \Gamma_O, \Gamma_{O^\perp}, \eta\rangle
$
such that
$\chi_3(O) = W_2$.
Thus we have
$\chi_3(O) = \psi \chi_2 \chi_1(O)$
so that
by Lemma 4.7,
$\chi_3^{-1} \psi \chi_2 \chi_1 \in \Gamma_O \Gamma_{O^\perp}$.
Since
$\chi_1, \chi_2,\chi_3\in 
\langle \Gamma_O, \Gamma_{O^\perp}, \eta\rangle,
$
we conclude that
$\psi \in 
\langle \Gamma_O, \Gamma_{O^\perp}, \eta\rangle.
$

Case 2: $O \cap \psi^{-1}(O^c) = \emptyset.$

Since $\psi(O)\cap O^c = \emptyset$,
there exists $\alpha \in \Gamma_A$ 
such that
$\alpha(\psi(O)) \subset O^c$.
Put
$V = \alpha(\psi(O))$.
By Lemma \ref{lem:3.18},
there exists 
$\chi 
\in 
\langle \Gamma_O, \Gamma_{O^\perp}, \eta\rangle
$
such that
$\chi(\psi(O)) = V$.
Put $\widetilde{\psi} = \chi\circ  \psi$.
As
$ \widetilde{\psi}(O) \subset O^c$,
by using Case 1 for 
$\widetilde{\psi}$ 
we see that
$\widetilde{\psi}$
and hence 
$\psi$ belongs to
$ 
\langle \Gamma_O, \Gamma_{O^\perp}, \eta\rangle.
$
\end{pf}
The following lemma is seen in \cite[Lemma 3.20]{GPS2}. 
\begin{lem}[{\cite[Lemma 3.20]{GPS2}}]\label{lem:3.20}
For a regular open set
$O \in RO(X_A)$, the following conditions are equivalent:
\begin{enumerate} 
\renewcommand{\labelenumi}{(\roman{enumi})}
\item 
$O$ is clopen.
\item 
For all $U \in RO(X_A)$ with $O \subset U$ and $O\ne U$,
we have $\overline{O} \subset U$.
\end{enumerate}
\end{lem}

\begin{lem}\label{lem:4.10}
Let $E \subset X_A$ be a nonempty closed of $X_A$
such that
$\gamma(E) \subset E$ for all $\gamma \in \Gamma_A$.
Then we have $E = X_A$.
This means that the action of $\Gamma_A$ on $X_A$ is minimal.
\end{lem}
\begin{pf}
Suppose that
$E \ne X_A$.
Take an arbitrary point $x \in E$.
Since $E^c$ is a nonempty open set, 
one may find clopen sets
$U, V \subset X_A$
such that
$$
x \in U, \qquad U \cap V = \emptyset, \qquad V \subset E^c.
$$
By Lemma 4.3,
there exists $\gamma \in \Gamma_A$ such that
$\gamma(U) \subset V$.
Since $\gamma(x) \in V \subset E^c$,
we have a contradiction, so that
$E = X_A$.
\end{pf}

A similar statement to the following lemma is seen in 
\cite[Lemma 3.21]{GPS2}.
The proof below is also similar to that of \cite[Lemma 3.21]{GPS2}.
We give the proof for the sake of completeness.
\begin{lem}[cf. {\cite[Lemma 3.21]{GPS2}}]\label{lem:3.21}
Let $O \subset X_A$ be a regular open set.
Then the following two conditions are equivalent:
\begin{enumerate} 
\renewcommand{\labelenumi}{(\roman{enumi})}
\item $ O$ is clopen. 
\item For a strong commuting pair $(H,K)$ of subgroups of $\Gamma_A$
such that $\Gamma_O \subset H$ with $\Gamma_O \ne H$,
the subgroup
$\langle H,\Gamma_{O^\perp} \rangle$ 
generated by 
$H$ and $\Gamma_{O^\perp}$ 
coincides with $\Gamma_A$.
\end{enumerate}
\end{lem}
\begin{pf}
(i) $\Longrightarrow$ (ii):
Assume that $O$ is clopen.
Let $(H,K)$ be a strong commuting pair of subgroups of $\Gamma_A$
satisfying $\Gamma_O \subset H$
with $\Gamma_A \ne H$.
Suppose that there exists
$\eta \in H$ such that
$\eta(O^c) \subset O$.
Then we have
$\eta^{-1}\Gamma_{O^\perp} \eta \subset \Gamma_O$
so that
$\Gamma_{O^\perp} \subset \eta \Gamma_O \eta^{-1} \subset H$.
Therefore by Lemma \ref{lem:3.9},
we have
$$
K  ={H }^\perp \subset (\Gamma_{O^\perp})^\perp = \Gamma_O \subset H.
$$
As $(H,K )$ is a strong commuting pair,
we have $K  = \{\id\}$ and hence $H = \Gamma_A$
so that
$\langle H,\Gamma_{O^\perp}\rangle = \Gamma_A$.
Suppose next that
$\eta(O^c) \cap O^c \ne \emptyset$
for all $\eta \in H$.
If  
$\eta(O) \subset O$ for all $\eta \in H$,
 then $\eta(O) = O$.
Then $\Gamma_O$ is a normal subgroup of $H$
so that
$(\Gamma_{O^\perp})^\perp =K$ and hence $H = \Gamma_O$.
Therefore we may assume that
there exists 
$
\eta \in H
$
such that
$$
\eta(O^c) \cap O^c \ne \emptyset, \qquad
\eta(O) \cap O^c \ne \emptyset.
$$
Hence by Lemma \ref{lem:3.19},
we have
$\langle \Gamma_O,\Gamma_{O^\perp}, \eta\rangle= \Gamma_A$
so that
$\langle H,\Gamma_{O^\perp}\rangle = \Gamma_A$.

(ii) $\Longrightarrow$ (i):
We will show that if $U \in RO(X_A)$
with $O \subset U$ and $O \ne U$,
then $\overline{O} \subset U$.
For $U \in RO(X_A)$,
consider the pair
$(\Gamma_U,\Gamma_{U^\perp})$
of subgroups
of $\Gamma_A$.
As $O \subset U$ with $O \ne U$, 
one has 
$\Gamma_O \subset \Gamma_U$
with
$\Gamma_O \ne \Gamma_U$
by Lemma \ref{lem:3.9} (i)
so that $\langle \Gamma_U,\Gamma_{O^\perp}\rangle = \Gamma_A$
by the condition (ii).
Assume that
$\overline{O} \cap U^c \ne \emptyset.$
Since
$\overline{O}$ is fixed by $\Gamma_{O^\perp}$
and
$U^c$ is fixed by $\Gamma_U$,
the closed set 
$\overline{O} \cap U^c $
is fixed by
 $\langle \Gamma_U,\Gamma_{O^\perp} \rangle$.
 Hence by the preceding lemma, 
 $\overline{O} \cap U^c  = X_A$
 a contradiction.
 Therefore we have
$\overline{O}\subset U$
and hence $O$ is clopen.
\end{pf}

Following \cite{GPS2},
we define condition (D3) for a strong commuting pair $(H,K)$
as follows:

\noindent
{\bf Definition.} 
A strong commuting pair $(H,K)$
of subgroups of $\Gamma_A$ is said to satisfy condition (D3)
if it satisfies the following two conditions (a) and (b):
\begin{enumerate} 
\renewcommand{\labelenumi}{(\alph{enumi})}
\item For a strong commuting pair $(H',K')$ of subgroups of $\Gamma_A$
such that $H \subset H'$ with $H \ne H'$,
the subgroup
$\langle H',K \rangle$ of $\Gamma_A$
generated by $H'$ and $K$ is equal to $\Gamma_A$.
\item For a strong commuting pair $(H^{\prime\prime},K^{\prime\prime})$ 
of subgroups of $\Gamma_A$
such that $K \subset K^{\prime\prime}$  
with $K \ne K^{\prime\prime}$,
the subgroup
$\langle H,K^{\prime\prime} \rangle$ of $\Gamma_A$
generated by $H$ and $K^{\prime\prime}$ 
is equal to $\Gamma_A$.
\end{enumerate}
We remark that 
in the statement of the above condition (b), 
the condition $K \subset K^{\prime\prime}$  with 
$K \ne K^{\prime\prime}$
is equivalent to
the condition
$H^{\prime\prime} \subset H$  with $H^{\prime\prime} \ne H$.
Therefore we have
\begin{prop}\label{prop:4.12}
Let $O \subset X_A$ be a regular open set. 
Then
$O$ is clopen if and only if 
the strong commuting pair 
$(\Gamma_O,\Gamma_{O^\perp})$
satisfies condition (D3).
\end{prop}
\begin{pf}
A regular open set 
$O \subset X_A$ is clopen if and only if 
$O^\perp$ is clopen.
Hence the assertion is clear by Lemma \ref{lem:3.21}.
\end{pf}

%%%%%%%%%%%%%%%%%%%%%%%%%%%%%%%%%%%%%%%%
%%%%%%%%%%%%%%%%%%%%%%%%%%%%%%%%%%
\section{Support of a strong commuting pair $(H,K)$}
%%%%%%%%%%%%%%%%%%%%%%%%%%%%%%%%%%%%%%%%%%

As in the preceding section, a clopen set $O$ (and $O^\perp$)
yields a strong commuting pair $(\Gamma_O,\Gamma_{O^\perp})$
of subgroups of $\Gamma_A$ satisfying condition (D3) (Proposition \ref{prop:4.12}).
In this section, we will conversely define a clopen set $P_H$ (and $P_K$)
of $X_A$ from a strong commuting pair $(H,K)$
satisfying condition (D3). 

Folllowing \cite[Definition 3.14]{GPS2},
we use the notations below.

For $\gamma \in \Gamma_A$,
denote by 
$X_A^\gamma$ the set of elements of $X_A$ fixed by $\gamma$:
\begin{equation*}
X_A^\gamma = \{ x \in X_A \mid \gamma(x) = x \}.
\end{equation*}
Denote by $P_\gamma$ the support of $\gamma$ defined by
\begin{equation*}
P_\gamma = \overline{(X_A^\gamma)^c}.
\end{equation*}
We note that 
$P_\gamma$ is a regular closed set and hence
$$
P_\gamma = \overline{(P_\gamma)^\circ} = (X_A^\gamma)^{\perp}.
$$
We in fact see that 
the inclusion relation
$\overline{(P_\gamma)^\circ} \subset P_\gamma$ is clear.
For the converse inclusion relation,
we see 
$
(\overline{(X_A^\gamma)^c})^\circ
\supset
 ((X_A^\gamma)^c)^\circ
=(X_A^\gamma)^c.
$
Hence 
$\overline{(P_\gamma)^\circ} \supset P_\gamma$.

For a subset $H \subset \Homeo(X_A)$,
define the support of $H$ as a closed subset of $X_A$ by
$$ 
P_H = \overline{\cup_{\eta \in H}(P_\eta)^\circ}.
$$

\noindent
{\bf Remark.}

{\bf 1.}
$P_H$ is a regular closed set such that
both $P_H$ and $(P_H)^\circ$ are $H^\perp$-invariant.

We in fact see that the inclusion relation
$\overline{(P_H)^\circ} \subset P_H$ is clear.  
Conversely, 
the inclusion relation
$\cup_{\eta \in H}(P_\eta)^\circ \subset P_H$
implies 
$\cup_{\eta \in H}(P_\eta)^\circ \subset (P_H)^\circ$
so that
$P_H \subset \overline{(P_H)^\circ}$.

We will next see that
$\xi(P_H) = P_H$ for $\xi \in H^\perp$. 
For $x \in X_A^\eta$ with $\eta \in H$, 
we have
$\eta(\xi(x)) = \xi(\eta(x)) = \xi(x)$
so that
$\xi(X_A^\eta) \subset X_A^\eta$
and similarly
$\xi^{-1}(X_A^\eta) \subset X_A^\eta$.
Hence we have
$\xi(X_A^\eta) = X_A^\eta$.
This implies
$\xi(P_\eta) = P_\eta$ for $\eta \in H$. 
As $\xi$ is a homeomorphism,
we have
$\xi((P_\eta)^\circ) = (P_\eta)^\circ$.
Hence we have
$\xi(P_H) = P_H$ for $\xi \in H^\perp$
and
$\xi((P_H)^\circ) = (P_H)^\circ$ for $\xi \in H^\perp$.

{\bf 2.}
For $\gamma \in \Gamma_A$, the set 
$P_\gamma$ is clopen. 

Its proof is the following.
As $\gamma$ is a cylinder map,
there exist $L\in {\Bbb N}$
and words
$\mu(\nu) = 
 (\mu_1(\nu),\dots,\mu_{k(\nu)}(\nu)) \in B_{k(\nu)}(X_A)
$
for
$\nu=(\nu_1,\dots,\nu_L)\in B_L(X_A)$
such that
%$$\gamma : X_A = \sqcup_{\nu \in B_L(X_A)}U_\nu 
%\longrightarrow X_A = \sqcup_{\mu(\nu) \in B_*(X_A)}U_{\mu(\nu)} 
%$$ and
$$
\gamma(\nu_1,\dots,\nu_L,x_{L+1},x_{L+2},x_{L+3},\dots )
=
(\mu_1(\nu),\dots,\mu_{k(\nu)}(\nu),x_{L+1},x_{L+2},x_{L+3},\dots )
$$
for $(\nu_1,\dots,\nu_L,x_{L+1},x_{L+2},x_{L+3},\dots ) \in U_{\nu}$.
Hence the set
$
(X_A^\gamma)^c
$ 
is a disjoint union of the following two clopen sets:
\begin{align*}
& \sqcup_{\nu \in B_L(X_A)} \{ U_\nu \mid L \ne k(\nu) \}, \\
& \sqcup_{\nu \in B_L(X_A)} \{ U_\nu \mid L = k(\nu),
(\nu_1,\dots,\nu_L) \ne (\mu_1(\nu),\dots,\mu_{k(\nu)}(\nu))\}.
\end{align*}
This implies that
$
(X_A^\gamma)^c
$ 
is clopen, 
and so is
$P_\gamma$.

{\bf 3.}
For $H \subset \Homeo(X_A)$ and 
$U \subset RO(X_A)$,
if $H \subset \Gamma_U$, then $P_H \subset \overline{U}$.

The proof is the following.
In general,
$H_1 \subset H_2 \subset \Homeo(X_A)$
implies
$P_{H_1} \subset P_{H_2}$.
Hence the condition
$H \subset \Gamma_U$
implies 
$P_H \subset P_{\Gamma_U}$.
By the lemma below,
one has 
$P_{\Gamma_U} = \overline{U}$
so that
$P_H \subset \overline{U}$.

\begin{lem}[cf. {\cite[Lemma 3.16]{GPS2}}]\label{lem:3.16} 
If $O\subset X_A$ 
is an open set, 
we have
$P_{\Gamma_O} = \overline{O}.$
Hence for a clopen set  
 $O\subset X_A$, 
we have
$P_{\Gamma_O} = O.$
\end{lem}
\begin{pf}
For $\eta \in \Gamma_O$
and
$x \in (X_A^\eta)^c$,
one sees that 
$\eta(x) \ne x$ so that 
$x \in O$.
Hence we have
$(X_A^\eta)^c \subset O$,
so that 
$P_\eta \subset \overline{O}$
and hence we have
$P_{\Gamma_O} \subset \overline{O}$.
For the converse inclusion relation,
take $x \in O$.
By Lemma \ref{lem:2.1}, 
there exist an open neighborhood $U$ of $x$
 and $\gamma \in \Gamma_U$ such that 
 $x \in U \subset O$ and
$ \gamma(x) \ne x$.
Take a clopen set $V$ such that
$x \in V \subset U$ and $\gamma(V) \cap V = \emptyset$.
We thus have
$V \subset (X_A^\gamma)^c$
so that 
$ V \subset \overline{(X_A^\gamma)^c} = P_\gamma$
and
$V \subset (P_\gamma)^\circ$.
Hence
$x \in \cup_{\gamma \in \Gamma_U}(P_\gamma)^\circ$.
Since
$\Gamma_U \subset \Gamma_O$,
we have
$
x \in 
\overline{\cup_{\gamma \in \Gamma_O}(P_\gamma)^\circ} = P_{\Gamma_O}$.
Therefore we have
$O \subset P_{\Gamma_O}$
and hence
$\overline{O} \subset P_{\Gamma_O}$.
 \end{pf}

We will next show that 
if $(H,K)$ satisfies condition (D3),
then the sets $P_H$ and $P_K$ are both clopen.
We provide a lemma.
\begin{lem}
Let $H$ be a subgroup of $\Gamma_A$.
\begin{enumerate} 
\renewcommand{\labelenumi}{(\roman{enumi})}
\item $\eta(y) = y$ for all $\eta \in H$ and $y \in (P_H)^c$.
\item Put $O = (P_H)^\circ$. Then we have 
$\zeta(O) = O$ and $\zeta(O^c) = O^c$ for all $\zeta \in H^\perp$.
\end{enumerate}
\end{lem}
\begin{pf}
(i)
Since we have
\begin{equation*}
(P_H)^c 
=
(\overline{\cup_{\eta \in H} (P_\eta)^\circ})^c
=
(\cap_{\eta \in H} ((P_\eta)^\circ)^c)^\circ
\subset
\cap_{\eta \in H} ((P_\eta)^\circ)^c
\end{equation*}
and
$
((P_\eta)^\circ)^c
=
\overline{(X_A^\eta)^\circ} \subset X_A^\eta,
$
we have
$
(P_H)^c 
\subset
\cap_{\eta \in H} X_A^\eta.
$

(ii)
We note that
$P_\eta$ is clopen for $\eta \in H$.
We have for $\zeta \in H^\perp$,
$\zeta(P_H) = \overline{\cup_{\eta \in H}\zeta(X_A^\eta)^c}$.
Since $\zeta$ commutes with $\eta \in H$,
we have
$
\zeta(X_A^\eta)
= X_A^\eta.
$
Hence
we see that
$\zeta(P_H) = P_H$ so that 
$\zeta(O) = O$ and $\zeta(O^c) = O^c$.
\end{pf}

\begin{lem}[cf. {\cite[Lemma 3.23]{GPS2}}]\label{lem:3.23}
Let $(H,K)$ be a strong commuting pair of subgroups of $\Gamma_A$
satisfying condition (D3).
Then the sets $P_H$ and $P_K$ are both clopen.
\end{lem}
\begin{pf}
Put
$O = (P_H)^\circ$.
As
$P_H$ is a regular closed set, 
$O$ is a regular open set satisfying
$\overline{O} = P_H.$
We have a strong commuting pair $(\Gamma_O,\Gamma_{O^\perp})$
of subgroups of $\Gamma_A$.
We  will prove that
$P_H = O$. 
By Lemma \ref{lem:3.16},
one knows that
$\overline{O} = P_{\Gamma_O}$ so that
$P_H = P_{\Gamma_O}.$
By the above lemma, 
$\eta(y) = y$ for all
$\eta \in H$ and $y \in (P_H)^c$.
As $O^c = \overline{(P_H)^c}$,
we have
$\eta(y) = y$
for all $\eta \in H$ and $y \in O^c$
so that   
$H \subset \Gamma_O$.
We have two cases.

Case 1: $H = \Gamma_O$.

Since $K = H^\perp = (\Gamma_O)^\perp = \Gamma_{O^\perp}$,
we have
$(H,K) = (\Gamma_O,\Gamma_{O^\perp})$. 
By the hypothesis, we see that
$(\Gamma_O,\Gamma_{O^\perp})$ satisfies condition (D3). 
Hence $O$ is clopen by Proposition \ref{prop:4.12}
so that $P_H = \overline{O} = O$ is clopen.

Case 2: $H\ne \Gamma_O$.

Suppose that $\overline{O} \cap O^c \ne \emptyset$.
Since $(H,K)$ satisfies condition (D3)
and
$H \subset\Gamma_O$
with $H \ne \Gamma_O$,
we have
$\langle \Gamma_O, K \rangle = \Gamma_A$.
As $O = (P_H)^\circ$, 
we know $\zeta(O) =O$ for all $\zeta \in H^\perp$
so that
$\zeta(\overline{O}) =\overline{O}$ for all $\zeta \in H^\perp =K$.
Hence the closed set 
$\overline{O} \cap O^c$
is globally invariant under 
$\langle \Gamma_O, K \rangle = \Gamma_A$.
By Lemma \ref{lem:4.10},
one knows that
$\overline{O} \cap O^c =\emptyset$
so that
$\overline{O} = O$.
Hence $O$ is clopen such that
$P_H = O$.

Symmetrically by considering $U = (P_K)^\circ$
one sees that 
$P_K$ is clopen.
\end{pf}

%%%%%%%%%%%%%%%%%%%%%%%%%%%%%%%%%%%%%%%%%%%%%%%%%%%%
%%%%%%%%%%%%%%%%%%%%%%%%%%%%%%%%%%%%%%%%%%%%%%%%%%%
\section{Dye Pairs}
%%%%%%%%%%%%%%%%%%%%%%%%%%%%%%%%%%%%%%%%%%%%%%%%%%
This section is devoted to proving Proposition \ref{prop:3.28} which asserts
that a strong commuting pair
$(H,K)$ with extra conditions (D4) and (D5) 
is of the form $(\Gamma_O, \Gamma_{O^\perp})$
for some clopen set $O$ of $X_A$.
The key lemma is Lemma \ref{lem:3.24} below.
We use the same conditions (D4) and (D5) as \cite[Definition 3.25]{GPS2}.

\noindent
{\bf Definition.} 
A strong commuting pair $(H,K)$ is said to satisfy conditions (D4) and (D5) if it satisfies the following  conditions:

(D4)
For a homeomorphism $\alpha \in \Gamma_A\backslash HK$, 
there exists 
       $\eta \in H$ with $\eta \ne \id$ 
(resp. $\kappa \in K$ with $\kappa \ne \id$)
such that
       $\alpha \eta \alpha^{-1} \in K$
(resp. $\alpha \kappa \alpha^{-1} \in H$).

(D5) If $N$ is a subgroup of $\Gamma_A$ with $N \ne \id$
such that $\eta N \eta^{-1} =N$ 
for all  $\eta \in H$, and $N \not \subset K$, 
(resp.    $\kappa N \kappa^{-1} =N$ 
for all  $\kappa \in K$, and $N \not \subset H$),
then $N \cap H \ne \{ \id \}$
(resp. $N \cap K \ne \{ \id \}$).

Following \cite{GPS2},
we will define the notion of Dye pair as in the following way.

\noindent
{\bf Definition.} 
A strong commuting pair $(H,K)$ satisfying condition (D3)
is said to be a {\it Dye pair}
if it satisfies the conditions (D4) and (D5).

The following two lemmas and their proofs
are similar to \cite[Lemma 3.26]{GPS2} and \cite[Lemma 3.27]{GPS2}.
We omit their proofs.

\begin{lem}[cf. {\cite[Lemma 3.26]{GPS2}}]\label{lem:3.26}
If $O$ is a clopen set,
then $(\Gamma_O, \Gamma_{O^\perp})$ is a Dye pair.
\end{lem}
%\begin{pf}
%It is enough to prove that 
% $(\Gamma_O, \Gamma_{O^\perp})$ satisfies conditions (D4) and (D5).
%Its proof is the same as the proof of \cite[Lemma3.26]{GPS2}. 
%\end{pf}
%The following lemma and its proof is the same as \cite[Lemma3.27]{GPS2}

\begin{lem}[cf. {\cite[Lemma 3.27]{GPS2}}]\label{lem:3.27} 
Let $(H,K)$ be a strong commuting pair of subgroups of $X_A$
satisfying conditions (D4) and (D5) such that 
$P_H = P_K = X_A$.
If $O \subset X_A$ is either an $H$- or $K$-invariant nonempty open set
of $X_A$, then $\overline{O} = X_A$.  
\end{lem}

For a subgroup
$H$ of $\Gamma_A$,
we define 
the continuous full group $[H]_c$
of $H$ as follows.
A homeomorphism $\gamma $ on $X_A$ belongs to
$[H]_c$ 
if there exist a finite clopen partition
$\sqcup_{i=1}^n U_i = X_A$ of $X_A$
and $\eta_i \in H$ such that
$X_A = \sqcup_{i=1}^n \eta_i(U_i)$
and
$\gamma(x) = \eta_i(x)$ for $x \in U_i$
(see \cite[Definition 2.2]{GPS2}).
Following \cite{GPS2},
for two homeomorphisms
$\alpha,\beta$ on $X_A$,
the closed set $F(\alpha,\beta)$  of $X_A$
is defined by
$$
F(\alpha,\beta) = \{ x \in X_A \mid \alpha(x) = \beta(x) \}.
$$

The statement of the following lemma is similar to 
\cite[Lemma 3.24]{GPS2}.
An invariant measure is used in the proof of \cite[Lemma 3.24]{GPS2}.
As in Proposition \ref{prop:2.4}, 
there is no $\Gamma_A$-invariant regular Borel probability measure on $X_A$,
so that we must modify its proof as in the following way.

\begin{lem}[cf. {\cite[Lemma 3.24]{GPS2}}] \label{lem:3.24}
Let $(H,K)$ be a strong commuting pair of subgroups of $X_A$
satisfying conditions (D4) and (D5).
Suppoe that 
$P_H = P_K = X_A$.
Then we have $[H]_c \cap [K]_c = \{ \id \}.$
\end{lem}
\begin{pf}
Suppose that
$[H]_c \cap [K]_c \ne \{ \id \}.$
Take $\alpha \in [H]_c \cap [K]_c$
with
$\alpha \ne \id$.
There exist
a nonempty clopen set
$U_0 \in CO(X_A)$
and homeomorphisms
$\eta_0 \in H, \kappa_0 \in K$ such that
$$
\alpha(x) = \eta_0(x) = \kappa_0(x), \qquad x \in U_0.
$$
Hence 
we have
$$
\emptyset \ne U_0 
\subset
F(\eta_0, \kappa_0)^\circ 
\subset F(\eta_0, \kappa_0) \ne X_A.
$$
If $\eta_1, \eta_2 \in H$
satisfy
$F(\eta_1,\eta_2)^\circ \ne \emptyset$,
we have for $x \in F(\eta_1,\eta_2)^\circ$ and $\kappa \in K$
$$
\eta_1(\kappa(x))=
\kappa(\eta_1(x)) =
\kappa(\eta_2(x)) = 
\eta_2(\kappa(x)),
$$
so that 
$F(\eta_1,\eta_2)^\circ$ is a $K$-invariant open set.
By Lemma \ref{lem:3.27},
we see that
$F(\eta_1,\eta_2)^\circ$ is dense in $X_A$
and hence 
$\eta_1 = \eta_2$ on $X_A$.
Let
$C(\eta_0)$ be the conjugacy class
$\{ \zeta \eta_0 \zeta^{-1} \in H \mid \zeta \in H\}$
of $\eta_0$ in $H$.
We similarly define the conjugacy class
$C(\kappa_0)$ of $\kappa_0$ in $K$.
For $\alpha, \beta \in C(\eta_0)$
with $\alpha \ne \beta$, 
we have 
$ F(\alpha,\beta)^\circ =\emptyset$.
For $\eta \in H$ and $\alpha\in C(\eta_0)$,
 we have
$$
\eta \alpha \eta^{-1}(\eta(x)) = \eta\alpha(x) =\eta(\kappa_0(x)) 
=\kappa_0(\eta(x)), \qquad  
x \in F(\alpha,\kappa_0)^\circ.
$$
Hence
$
\eta(F(\alpha,\kappa_0)^\circ) 
\subset 
F(\eta \alpha\eta^{-1},\kappa_0)^\circ
$
and 
$
\eta^{-1}(F(\eta\alpha\eta^{-1},\kappa_0)^\circ) 
\subset F(\alpha,\kappa_0)^\circ.
$
We thus have with symmetric discussions for $C(\kappa_0)$,
\begin{enumerate} 
\renewcommand{\labelenumi}{(\roman{enumi})}
\item 
\ \ \ $F(\alpha,\kappa_0)^\circ \cap F(\beta,\kappa_0)^\circ
\subset F(\alpha,\beta)^\circ =\emptyset$ for $\alpha, \beta \in C(\eta_0)$ with $\alpha \ne \beta$,  

$F(\eta_0,\alpha)^\circ \cap F(\eta_0,\beta)^\circ
\subset F(\alpha,\beta)^\circ =\emptyset$ for $\alpha, \beta \in C(\kappa_0)$ with $\alpha \ne \beta$.
\item
\ \ \ $\eta(F(\alpha,\kappa_0)^\circ) =
 F(\eta \alpha \eta^{-1},\kappa_0)^\circ
 $ for $\eta \in H, \alpha \in C(\eta_0)$, 

$\kappa(F(\eta_0,\beta)^\circ) =
 F(\eta_0,\kappa\beta \kappa^{-1})^\circ
 $ for $\kappa \in K, \beta \in C(\kappa_0)$.
\end{enumerate}
Let
\begin{equation*}
N_H = \langle \Gamma_{F(\alpha,\kappa_0)^\circ}:\alpha \in C(\eta_0) \rangle 
\quad
\text{ and }
\quad
N_K = \langle \Gamma_{F(\eta_0,\beta)^\circ}: \beta \in C(\kappa_0) \rangle
\end{equation*}
be the subgroup of $\Gamma_A$ generated by elements of  
$\cup_{\alpha \in C(\eta_0)}\Gamma_{F(\alpha,\kappa_0)^\circ}$
and
that of  $\Gamma_A$ generated by elements of  
$\cup_{\beta \in C(\kappa_0)}\Gamma_{F(\eta_0,\beta)^\circ}$
respectively.
Since for $\eta\in H, \, \alpha \in C(\eta_0)$
\begin{equation}
\eta \Gamma_{F(\alpha,\kappa_0)^\circ} \eta^{-1} 
=
\Gamma_{\eta(F(\alpha,\kappa_0)^\circ)}
= 
\Gamma_{F(\eta\alpha\eta^{-1},\kappa_0)^\circ}, \label{eqn:6.1}
\end{equation}
we know that
\begin{equation*}
\eta N_H \eta^{-1} = N_H
\text{ for all } 
\eta \in H,
\text{ and similarly }
\kappa N_K \kappa^{-1} = N_K
\text{ for all } 
\kappa \in K.
\end{equation*}
We will first show that
$N_K \not\subset H$ and $N_H \not\subset K$.
Suppose that
$N_K \subset H$.
We then have 
$\Gamma_{F(\eta_0,\beta)^\circ} \subset H$
for all $\beta \in C(\kappa_0)$,
in particular,
$\Gamma_{F(\eta_0,\kappa_0)^\circ} \subset H.$
By \eqref{eqn:6.1},
we see that
$
\Gamma_{F(\alpha,\kappa_0)^\circ}  \subset H$
for $\alpha \in C(\eta_0)$.
 Hence we have
$N_H \subset H$.
Put for $x \in X_A$
\begin{align*}
\widetilde{\eta}_0(x) 
& = 
{
\begin{cases}
\eta_0(x) & \text{ if } x \in F(\eta_0,\kappa_0)^\circ, \\
        x & \text{ otherwise,}
\end{cases}
} \\ 
\widetilde{\kappa}_0(x) 
& = 
{
\begin{cases}
\kappa_0(x) & \text{ if } x \in F(\eta_0,\kappa_0)^\circ, \\
          x & \text{ otherwise.}
\end{cases}
}
\end{align*}
Since
$\eta_0(F(\eta_0,\kappa_0)^\circ) =F(\eta_0,\kappa_0)^\circ$
and
$\widetilde{\eta}_0|_{(F(\eta_0,\kappa_0)^\circ)^c} =\id$,
one sees that
$\widetilde{\eta}_0 \in \Gamma_{F(\eta_0,\kappa_0)^\circ} \subset H$.
We will next show 
$\widetilde{\kappa}_0 \in K$.
Take an arbitrary $\eta \in \Gamma_{F(\eta_0,\kappa_0)^\circ}$.
For
$x \in F(\eta_0,\kappa_0)^\circ$,
we have
\begin{equation*}
\eta \widetilde{\kappa}_0(x)
=\eta \kappa_0(x)
=\kappa_0\eta(x)
=\widetilde{\kappa}_0\eta(x).
\end{equation*}  
For
$x \in (F(\eta_0,\kappa_0)^\circ)^c$,
we have
\begin{equation*}
\eta \widetilde{\kappa}_0(x)
=\eta (x)
=\widetilde{\kappa}_0\eta(x).
\end{equation*}  
Hence we have
%$ \eta \widetilde{\kappa}_0(x)=\widetilde{\kappa}_0\eta(x)$ 
%for all $x \in X_A$ so that
\begin{equation} 
\eta \widetilde{\kappa}_0
=\widetilde{\kappa}_0\eta
\text{ for all }
\eta \in \Gamma_{F(\eta_0,\kappa_0)^\circ}. \label{eqn:etakappa1}
\end{equation}
Take an arbitrary $\eta \in \Gamma_{F(\alpha,\kappa_0)^\circ}$
with $\alpha \in C(\eta_0)$ 
such that
$
F(\alpha,\kappa_0)^\circ \cap F(\eta_0,\kappa_0)^\circ = \emptyset.
$
For
$x \in F(\eta_0,\kappa_0)^\circ$,
we have
\begin{equation*}
\eta \widetilde{\kappa}_0(x)
=\eta \kappa_0(x)
=\kappa_0\eta(x)
=\kappa_0(x)
=\widetilde{\kappa}_0(x)
=\widetilde{\kappa}_0\eta(x).
\end{equation*}  
For
$x \in (F(\eta_0,\kappa_0)^\circ)^c$,
we have
\begin{equation*}
\eta \widetilde{\kappa}_0(x)
=\eta (x)
=\widetilde{\kappa}_0\eta(x).
\end{equation*}  
Hence we have
%$ \eta \widetilde{\kappa}_0(x)=\widetilde{\kappa}_0\eta(x)$ 
%for all $x \in X_A$ so that
\begin{equation} 
\eta \widetilde{\kappa}_0
=\widetilde{\kappa}_0\eta
\text{ for all }
\eta \in \Gamma_{F(\alpha,\kappa_0)^\circ}
\text{ with }
F(\alpha,\kappa_0)^\circ \cap F(\eta_0,\kappa_0)^\circ = \emptyset.
 \label{eqn:etakappa2}
\end{equation}
By \eqref{eqn:etakappa1} and \eqref{eqn:etakappa2},
we have
$\widetilde{\kappa}_0 \in (N_H)^\perp =K$
so that
$$
\id \ne \widetilde{\eta}_0 
= \widetilde{\kappa}_0 
\in H \cap K =\{ \id \}
$$
a contradiction.
Therefore we conclude that
$N_K \not\subset H$ and similarly $N_H \not\subset K$.

Now $N_K \ne \{ \id \}$
such that
$\kappa N_K \kappa^{-1} = N_K$ for all $\kappa \in K$
and
$N_K \not\subset H$.
By condition (D5), we see that
$N_K \cap K \ne \{ \id \}$
and similarly
$N_H \cap H \ne \{ \id \}$.
We set 
$$
\widetilde{N}_K = N_K \cap K, \qquad 
\widetilde{N}_H = N_H \cap H
$$
so that 
$\widetilde{N}_K$ is a normal subgroup of $K$
and
$\widetilde{N}_H$ is a normal subgroup of $H$.
By condition (D2), we see that
\begin{equation*}
(\widetilde{N}_K)^\perp = H, \qquad
(\widetilde{N}_H)^\perp = K.
\end{equation*}
Let
$\widetilde{\eta}_0$ and
$\widetilde{\kappa}_0$
be previously defined homeomorphisms on $X_A$.
In this setting we will show that
$\widetilde{\eta}_0 \in H$ and
$\widetilde{\kappa}_0 \in K$ as follows:

Since $F(\eta_0,\beta)^\circ, F(\eta_0,\beta')^\circ$
for $\beta, \beta' \in C(\kappa_0)$
are disjoint or equal,
we have for any 
$\gamma \in \Gamma_{F(\eta_0, \alpha)^0}, \alpha \in C(\kappa_0)$
\begin{equation*}
\gamma(F(\eta_0,\beta)^\circ) = F(\eta_0,\beta)^\circ
\quad
\text{ for }
\quad
 \beta \in C(\kappa_0).
\end{equation*}
Hence each 
$F(\eta_0,\beta)^\circ$
for $\beta\in C(\kappa_0)$
is globally invariant under
$N_K$ and $\widetilde{N}_K$.
Let $\kappa \in \widetilde{N}_K$ be an arbitrary element.
For
$x \in F(\eta_0,\kappa_0)^\circ$,
we have
\begin{equation*}
\kappa \widetilde{\eta}_0(x)
=\kappa\eta_0(x)
=\eta_0\kappa(x)
=\widetilde{\eta}_0\kappa(x).
\end{equation*}  
For
$x \in (F(\eta_0,\kappa_0)^\circ)^c$,
we have
\begin{equation*}
\kappa \widetilde{\eta}_0(x)
=\kappa (x)
=\widetilde{\eta}_0\kappa(x).
\end{equation*}  
Hence we have
$ 
\kappa \widetilde{\eta}_0
=\widetilde{\eta}_0\kappa
$
for all
$
\kappa \in \widetilde{N}_K
$
so that
we have
$\widetilde{\eta}_0 \in (\widetilde{N}_K)^\perp =H$,
symmetrically 
$\widetilde{\kappa}_0 \in (\widetilde{N}_H)^\perp =K$.
Hence we have
$$
\id \ne \widetilde{\eta}_0 
= \widetilde{\kappa}_0 
\in H \cap K =\{ \id \}
$$
a contradiction.
Therefore we conclude that
$[H]_c \cap [K]_c =\{ \id \}$.
\end{pf}

\begin{lem} \label{lem:3.24'}
Let $(H,K)$ be a strong commuting pair of subgroups of $X_A$
satisfying conditions (D4) and (D5).
Suppose that 
$P_H = P_K = X_A$.
Then there exists $\alpha \in [H]_c$
such that the fixed point set 
$X_A^\alpha$ is not $K$-invariant.
\end{lem}
\begin{pf}
As $P_H = X_A$, 
there exists $\gamma \in H$ and
$x \in X_A$ such that
$\gamma(x) \ne x$.
One may take a clopen neighborhood $V$ of $x$ such that
\begin{equation}
\gamma(V) \cap V =\emptyset, \qquad
\overline{\gamma(V) \cup V} \ne X_A. \label{eqn:gamma}
\end{equation}
Define $\alpha \in [H]_c$
by setting for $x \in X_A$
\begin{equation*}
\alpha(x) =
\begin{cases}
     \gamma(x) & \text{ for } x\in V,\\
\gamma^{-1}(x) & \text{ for } x\in \gamma(V),\\
             x & \text{ for } x\in (V \cup \gamma(V))^c
\end{cases}
\end{equation*}
so that
$(X_A^\alpha)^c = \gamma(V) \cup V$.
If $X_A^\alpha$ is $K$-invariant,
$(X_A^\alpha)^c$ is a $K$-invariant nonempty open set of $X_A$.
By Lemma \ref{lem:3.27}, 
one sees that
$(X_A^\alpha)^c$ is dense in $X_A$.
This contradicts to \eqref{eqn:gamma}.
Hence 
$X_A^\alpha$ is not $K$-invariant.
 \end{pf}

Before reaching the final proposition, we provide a lemma below.
\begin{lem}[cf. {\cite[Lemma 3.17]{GPS2}}]\label{lem:3.17}
Let
$(H,K)$ be a strong commuting pair of $\Gamma_A$.
\begin{enumerate} 
\renewcommand{\labelenumi}{(\roman{enumi})}
\item 
For a nonempty clopen set $U\subset P_H$ such that
$\gamma(U) =U $ for all $\gamma \in \langle H,K\rangle$, 
then $U = P_H$.
\item 
For a nonempty clopen set $V\subset P_K$ such that
$\gamma(V) =V $ for all $\gamma \in \langle H,K\rangle$, 
then $V = P_K$.
\end{enumerate}
\end{lem}
\begin{pf}
The proof is the same as that of \cite[Lemma 3.17]{GPS2}.
\end{pf}

Therefore we have
\begin{prop}[cf. {\cite[Proposition 3.28]{GPS2}}]\label{prop:3.28} 
Let $(H,K)$ be a Dye pair of subgroups of $X_A$.
Then 
$P_H$ and $P_K$ are clopen sets of $X_A$ such that 
\begin{equation*}
P_H = (P_K)^\perp, \qquad 
P_K = (P_H)^\perp
\quad
\text{ and }
\quad
(H,K) = (\Gamma_{P_H}, \Gamma_{P_K}).
\end{equation*}
\end{prop}
\begin{pf}
By Lemma \ref{lem:3.23}, 
both $P_H$ and $P_K$ are clopen.
Suppose that
$P_H \cap P_K \ne \emptyset$.
Since 
$\eta(\cap_{\zeta \in H}X_A^\zeta) 
=\cap_{\zeta \in H}X_A^\zeta$
for $\eta \in H$
and
hence
$\eta((\cap_{\zeta \in H}X_A^\zeta)^c) =(\cap_{\zeta \in H}X_A^\zeta)^c$,
we have
$\eta(P_H) = P_H$ for $\eta \in H$.
Similarly 
$\kappa(P_K) = P_K$ for $\kappa \in K$.
As
$P_H$ is $H^\perp(=K)$-invariant
and $P_K$ is $K^\perp(=H)$-invariant,
the set
$P_H \cap P_K$ is an $H$-and $K$-invariant clopen set.
By Lemma \ref{lem:3.17},
we have
$P_H \cap P_K = P_H = P_K$.
As in the proof of Lemma \ref{lem:3.23},
the inclusion relations
$H \subset \Gamma_{P_H},
K \subset \Gamma_{P_K}
$
hold.
Hence we have
\begin{equation*}
\Gamma_{(P_H)^c}
=\Gamma_{(P_H)^\perp}
=(\Gamma_{P_H})^\perp
\subset H^\perp =K 
\subset \Gamma_{P_K}
=\Gamma_{P_H}
\end{equation*}
so that
$(P_H)^c = \emptyset$.
Therefore we have
$P_H = P_K = X_A$.
By Lemma \ref{lem:3.24}, 
we have
\begin{equation}
[H]_c \cap [K]_c = \{ \id \}. \label{eqn:HK}
\end{equation}
By Lemma \ref{lem:3.24'},
there exists
$\alpha \in [H]_c$ 
such that
the fixed point set
$X_A^\alpha$ is not $K$-invariant.
This implies $\alpha \not\in H$.
If 
$\alpha = \eta \kappa $ for some $\eta \in H$ 
and $\kappa \in K$,
we see that 
$\eta^{-1} \alpha =\kappa $ belongs to $[H]_c \cap [K]_c$.
Hence
\eqref{eqn:HK}
implies 
$\eta^{-1}\alpha = \kappa = \id$
and $\alpha = \eta \in H$, a contradiction,
so that 
$\alpha$ does not belong to
the subgroup
$HK$.
By condition (D4), one may find 
$\eta \in H \backslash \{\id \}$
such that
$\alpha \eta \alpha^{-1} \in K$.
Since
$\alpha \eta \alpha^{-1} \in [H]_c \cap K 
\subset [H]_c \cap [K]_c =\{\id\},
$
we obtain $ \eta = \id$ a contradiction.
Thus we conclude 
$$
P_H \cap P_K= \emptyset.
$$
This implies
$P_H \subset (P_K)^c =(P_K)^\perp$.
It then follows that
$$
\Gamma_{P_H} \subset \Gamma_{(P_K)^\perp}
=(\Gamma_{P_K})^\perp
\subset K^\perp =H 
\subset \Gamma_{P_H}.
$$
Therefore we have 
$
\Gamma_{P_H}  =H,
P_H = (P_K)^\perp
$
 and similarly 
$
\Gamma_{P_K}  =K,
P_K = (P_H)^\perp.
$
\end{pf}

%%%%%%%%%%%%%%%%%%%%%%%
%%%%%%%%%%%%%%%%%%%%%%%%%%%%%%%%%%%%%%%
\section{Main result and its corollaries}
%%%%%%%%%%%%%%%%%%%%%%%%%%%%%%%%%%%%%%%%%%%%%%%
Let $A,B$ be two irreducible square matrices 
with entries in $\{0,1\}$ satisfying condition (I).
\begin{lem}\label{lem:7.1}
Suppose that there exists an isomorphism
$\alpha: \Gamma_A \longrightarrow \Gamma_B$
of groups.
Then for two subgroups $H,K $ of $\Gamma_A$, 
the following two conditions are equivalent:
\begin{enumerate} 
\renewcommand{\labelenumi}{(\roman{enumi})}
\item 
$(H,K)$ is a Dye pair of $\Gamma_A$.
\item
$(\alpha(H),\alpha(K))$ is a Dye pair of $\Gamma_B$.
\end{enumerate}
\end{lem}
\begin{pf}
We will show the implication 
(i) $\Longrightarrow$ (ii).
Let 
$(H,K)$ be a Dye pair of $\Gamma_A$.
Since 
$\alpha:\Gamma_A \longrightarrow \Gamma_B$
is an isomorphism of groups, 
one knows that
$(\alpha(H),\alpha(K))$ satisfies (D1), (D2) and (D3),
and hence it is a strong commuting pair of $\Gamma_B$
satisfying condition (D3). 
The conditions (D4) and (D5) are also determined by group structure,
so that $(\alpha(H),\alpha(K))$ satisfies (D4) and (D5).
Hence
$(\alpha(H),\alpha(K))$ is a Dye pair of $\Gamma_B$.
\end{pf}
An isomorphism $\alpha:\Gamma_A\longrightarrow \Gamma_B$
of groups is said to be {\it spatial} 
if there exists a homeomorphism
$h: X_A \longrightarrow X_B$ 
 such that $\alpha(\gamma) = h \circ \gamma \circ h^{-1}$
 for $\gamma \in \Gamma_A$.
%\begin{equation}
%h \circ \Gamma_A \circ h^{-1} = \Gamma_B. \label{eqn:spatial}
%\end{equation}
We arrive at the main result of the paper.

\begin{thm}
Let $A,B$ be two irreducible square matrices 
with entries in $\{0,1\}$ satisfying condition (I).
Then every group isomorphism 
$\alpha: \Gamma_A \longrightarrow \Gamma_B$
is spatial.
\end{thm}
\begin{pf}
The proof is achieved by constructing
a Boolean isomorphism $\varphi : CO(X_A) \longrightarrow CO(X_B)$
satisfying
\begin{equation*}
\alpha(\gamma) \circ \varphi = \varphi\circ \gamma, 
\qquad \gamma \in \Gamma_A. 
\end{equation*} 
The construction of $\varphi$ is due to 
Lemma \ref{lem:3.26},
Propositions \ref{prop:3.28} and 
Lemma \ref{lem:7.1}.
The detailed proof  is the same as the proof of \cite[Theorem 4.2]{GPS2}.
\end{pf}

Let us denote by $\OA$ 
the Cuntz-Krieger algebra for the matrix $A$,
and $\DA$ its canonical maximal abelian subalgebra of $\OA$ (\cite{CK}).

\begin{cor}
Let $A,B$ be two irreducible square matrices 
with entries in $\{0,1\}$ satisfying condition (I).
Then the following three conditions are equivalent:
\begin{enumerate} 
\renewcommand{\labelenumi}{(\roman{enumi})}
\item
The groups $\Gamma_A$ and $\Gamma_B$ 
 are isomorphic as abstract groups.
%\item The action of $\Gamma_A$ on $X_A$ and $\Gamma_B$ on $X_B$ are spatial.
\item
The one-sided topological Markov shifts
$(X_A,\sigma_A)$ and $(X_B,\sigma_B)$  
are continuously orbit equivalent.
\item
There exists an isomorphism $\Psi: \OA \rightarrow  \OB$
 such that $\Psi(\DA) = \DB$.
\end{enumerate}
\end{cor}

Let $N$ and $M$ be the size of matrix $A$ and that of $B$ respectively.
Denote by $I_N$ and by $I_M$ 
the identity matrix of size $N$ and that of size $M$ respectively.  
In \cite{MaPAMS},
under the condition that
$\det(A-I_N)\det(B-I_M) \ge 0$,
an isomorphism between Cuntz-Krieger algebras 
induces an isomorphism between them which preserves their canonical 
maximal abelian subalgebras.
Hence we have
\begin{cor}
Let $A,B$ be two irreducible square matrices 
with entries in $\{0,1\}$ satisfying condition (I).
Suppose that
$\det(A-I_N)\det(B-I_M) \ge 0$.
Then the following two conditions are equivalent:
\begin{enumerate} 
\renewcommand{\labelenumi}{(\roman{enumi})}
\item
The groups $\Gamma_A$ and $\Gamma_B$ 
 are isomorphic as abstract groups.
\item
The Cuntz-Krieger algebras
$\OA$ and $\OB$ are isomorphic.
\end{enumerate}
\end{cor}
By using classification theorem for Cuntz-Krieger algebras 
obtained by M. R{\o}rdam
\cite{Ro}(cf. \cite{Ro2}),
one may classify the continuous full groups 
in terms of the underlying matrices under the determinant contition
$\det(A-I_N) \det(B-I_M) \ge 0$ as follows:
\begin{cor}
Let $A,B$ be two irreducible square matrices 
with entries in $\{0,1\}$ satisfying condition (I).
Suppose that
$\det(A-I_N)\det(B-I_M) \ge 0$.
The groups $\Gamma_A$ and $\Gamma_B$ 
 are isomorphic as abstract groups
 if and only if
 there exists an isomorphism
$\varPhi: {\Bbb Z}^N / (A^t -I_N) {\Bbb Z}^N
 \longrightarrow
  {\Bbb Z}^M / (B^t -I_M) {\Bbb Z}^M
  $
  such that
$\varPhi([1,\dots,1]) = [1,\dots,1].$
\end{cor}
Therefore we know that thre are many mutually nonisomorphic 
continuous full groups of one-sided topological Markov shifts
such as the following corollary. 
\begin{cor}
Let $N,M$ be positive integers such that $N,M > 1$.
Denote by 
$\Gamma_{[N]}$ and $ \Gamma_{[M]}$ the continuous full groups of the one-sided full $N$-shift and $M$-shift respectively.  
Then 
$\Gamma_{[N]}$ and $ \Gamma_{[M]}$ are isomorphic as abstract groups if and only if
$N = M$.
\end{cor}

{\it Acknowledgments:}
This work was supported by JSPS Grant-in-Aid for Scientific Reserch 
((C), No 23540237).

\end{document}